\documentclass[11pt]{article}
\usepackage{amsmath}
\usepackage{amssymb}
\usepackage{amsthm}
\usepackage{latexsym}
\usepackage{graphicx}
\usepackage{hyperref}

\newtheorem{thm}{Theorem}[section]
\newtheorem{cor}[thm]{Corollary}
\newtheorem{lem}[thm]{Lemma}
\newtheorem{prop}[thm]{Proposition}

\newcommand{\mb}{\mathbf}
\newcommand{\mh}{\mathbb}
\newcommand{\mr}{\mathrm}
\newcommand{\mc}{\mathcal}
\newcommand{\mf}{\mathfrak}
\newcommand{\ds}{\displaystyle}
\newcommand{\scs}{\scriptstyle}
\newcommand{\ts}{\textstyle}

\newcommand{\es}{\emptyset}
\newcommand{\af}{\mr{aff}}

\newcommand{\inp}[2]{\langle #1 \,,\, #2 \rangle}
\newcommand{\prefix}[3]{\vphantom{#3}#1#2#3}

\begin{document}
\
\vspace{2cm}
\begin{center}
\LARGE \textbf{Parabolically induced representations\\ of graded Hecke algebras}\\[5mm]
Maarten Solleveld \\[2cm]
\normalsize Mathematisches Institut,
Georg-August-Universit\"at G\"ottingen\\
Bunsenstra\ss e 3-5, 37073 G\"ottingen, Germany\\
email: Maarten.Solleveld@mathematik.uni-goettingen.de \\
October 2009\\[3cm]
\end{center}
\textbf{Abstract.}\\
We study the representation theory of graded Hecke algebras, starting from scratch and 
focusing on representations that are obtained by induction from a discrete series representation 
of a parabolic subalgebra. We determine all intertwining operators between such parabolically 
induced representations, and use them to parametrize the irreducible representations.
Finally we describe the spectrum of a graded Hecke algebra as a topological space.\\[2mm]
\textbf{Mathematics Subject Classification (2000).} \\
20C08\\

\newpage

\tableofcontents

\vspace{4mm}

\section*{Introduction}
\addcontentsline{toc}{section}{Introduction}

This article aims to describe an essential part of the representation theory of
graded Hecke algebras. We will do this in the spirit of Harish-Chandra, 
with tempered representations and parabolic subalgebras, and we ultimately try to 
understand the spectrum from the noncommutative geometric point of view.

Graded Hecke algebras were introduced by Lusztig to facilitate the study of representations 
of affine Hecke algebras and simple groups over $p$-adic fields \cite{Lus-Gr,Lus1}.
While the structure of Hecke algebras of simple $p$-adic groups is not completely understood,
graded Hecke algebras have a very concrete definition in terms of generators and relations. 
According to Lusztig's reduction theorem, the relation between an affine Hecke algebra and
its graded version is similar to that between a Lie group and its Lie algebra, the graded
Hecke algebra is a kind of linearization of the affine one. Hence it is somewhat easier to
study, while still containing an important part of the representation theory. In \cite{Lus2}
Lusztig succeeded in parametrizing its irreducible representations in terms of ``cuspidal
local data''. Unfortunately it is unclear to the author how explicit this is, or can be made.

Alternatively a graded Hecke algebra $\mh H$ can be regarded as a deformation of 
$S (\mf t^* ) \rtimes W$, where $\mf t^*$ is some complex vector space containing a 
crystallographic root system $R$ with Weyl group $W$. The algebra $\mh H$ contains
copies of $S(\mf t^*)$ and of the group algebra $\mh C [W]$, and the multiplication 
between these parts depends on several deformation parameters $k_\alpha \in \mh C$. 

This point of view shows that there is a relation between the representation theory of 
$\mh H$ and that of $S (\mf t^* ) \rtimes W$. Since the latter is completely described by 
Clifford theory, this could offer some new insights about $\mh H$.
Ideally speaking, one would like to find data that parametrize both the irreducible 
representations of $\mh H$ and of $S (\mf t^* ) \rtimes W$. This might be asking too much, 
but it is possible if we are a little more flexible, and allow virtual representations. 

For this purpose it is important to understand the geometry of the spectrum of $\mh H$. 
(By the spectrum of an algebra $A$ we mean the set Irr$(A)$ of equivalence classes of 
irreducible representations, endowed with the Jacobson topology.) So we want a 
parametrization of the spectrum which highlights the geometric structure. Our approach 
is based on the discrete series and tempered representations of $\mh H$.
These are not really discrete or tempered in any obvious sense, the terminology is merely 
inspired by related classes of representations of affine Hecke algebras and of reductive groups
over $p$-adic fields. These $\mh H$-representations behave best when all deformation 
parameters are real, so will assume that in the introduction. 

Our induction data are triples $(P,\delta ,\lambda )$ consisting of a set of simple roots $P$, 
a discrete series representation $\delta$ of the parabolic subalgebra corresponding to $P$, 
and a character $\lambda \in \mf t$. From this we construct a ``parabolically induced 
representation'' $\pi (P,\delta ,\lambda )$. Every irreducible $\mh H$-representation appears as
a quotient of such a representation, often for several induction data 
$(P,\delta ,\lambda )$.

Usually $\pi (P,\delta ,\lambda )$ is reducible, and more specifically it is tempered and 
completely reducible if the real part $\Re (\lambda )$ of $\lambda$ is 0. There is a Langlands 
classification for graded Hecke algebras, analogous to the original one for reductive groups, 
which reduces the classification problem to irreducible tempered $\mh H$-representations. 
Thus three tasks remain:
\begin{description}
\item[a)] Find all equivalence classes of discrete series representations.
\item[b)] Determine when $\pi (P,\delta ,\lambda)$ and $\pi (Q,\sigma ,\mu )$ have common
irreducible constituents.
\item[c)] Decompose $\pi (P,\delta ,\lambda )$ into irreducibles, at least when $\Re (\lambda ) = 0$.
\end{description}
In \cite{OpSo} a) is solved, via affine Hecke algebras. 
For b) and c) we have to study $\mr{Hom}_{\mh H}(\pi (P,\delta ,\lambda ) , 
\pi (Q,\sigma ,\mu ))$, which is the main subject of the present paper. In analogy with 
Harish-Chandra's results for reductive groups, intertwiners between parabolically induced 
$\mh H$-modules should come from suitable elements of the Weyl group $W$. Indeed $W$ acts 
naturally on our induction data, and for every $w \in W$ such that $w(P)$ consists of simple roots 
we construct an intertwiner
\[
\pi (w,P,\delta ,\lambda ) : \pi (P,\delta ,\lambda ) \to \pi (w(P,\delta ,\lambda )) \,,
\]
which is rational as a function of $\lambda$. In Theorem \ref{thm:7.1} we prove

\begin{thm}\label{thm:0.1}
For $\Re (\lambda ) = 0$ the operators
\[
\{ \pi (w,P,\delta ,\lambda ) : w (P,\delta ,\lambda ) \cong (Q,\sigma ,\mu ) \}
\]
are regular and invertible, and they span 
$\mr{Hom}_{\mh H}(\pi (P,\delta ,\lambda) ,\pi (Q,\sigma ,\mu ) )$.
\end{thm}

The proof runs via affine Hecke algebras and topological completions thereof, making
this a rather deep result. In particular this is one of the points where the parameters
$k_\alpha$ have to be real.

Thus the answer to b) is that the packets of irreducible constituents of two tempered
parabolically induced representations are either disjoint or equivalent. Moreover 
this is already detected by $W$. However, it is hard to determine these intertwining 
operators explicitly, and sometimes they are linearly dependent, so this solves 
only a part of c). 

Furthermore it is possible to incorporate the Langlands classification in this picture.
Therefore we relax our condition on $\lambda$, requiring only that $\Re (\lambda )$
is contained in a certain positive cone. In Proposition \ref{prop:7.3}.c 
we generalize the above:

\begin{thm}\label{thm:0.2}
Suppose that $(P,\delta ,\lambda)$ and $(Q,\sigma ,\mu )$ are induction data, with 
$\Re (\lambda )$ and $\Re (\mu )$ positive. The irreducible quotients of $\pi 
(P,\delta ,\lambda )$ for such induction data exhaust the spectrum of $\mh H$. The operators 
\[
\{ \pi (w,P,\delta ,\lambda ) : w (P,\delta ,\lambda ) \cong (Q,\sigma ,\mu ) \}
\]
are regular and invertible, and they span 
$\mr{Hom}_{\mh H}(\pi (P,\delta ,\lambda ) ,\pi (Q,\sigma ,\mu ))$.
\end{thm}

This is about as far as the author can come with representation theoretic methods. 
To learn more about the number of distinct irreducibles contained in a parabolically 
induced representation, we call on noncommutative geometry. The idea here is not to 
consider just a finite packet of $\mh H$-representations, but rather to study its spectrum 
Irr$(\mh H)$ as a topological space. Since $\mh H$ is a deformation of $S(\mf t^*) \rtimes W$,
there should be a relation between the spectra of these two algebras. An appropriate
theory to make this precise is periodic cyclic homology $HP_*$, since there is an
isomorphism $HP_* (\mh H) \cong HP_* (S(\mf t^*) \rtimes W)$ \cite{Sol2}. As the 
periodic cyclic homology of an algebra can be regarded as a kind of cohomology of 
the spectrum of that algebra, we would like to understand the geometry of Irr$(\mh H)$ 
better. The center of $\mh H$ is $S(\mf t^*)^W = \mh C [\mf t / W]$ \cite{Lus-Gr},
so Irr$(\mh H)$ is in first approximation, the vector space $\mf t$
modulo the finite group $W$. A special case of Theorem \ref{thm:9.1} tells us that

\begin{thm}\label{thm:0.3}
Let $\pi (P,\delta ,\lambda ) = V_1 \oplus V_2$ be a decomposition of $\mh H$-modules, 
and suppose that $\lambda'$ satisfies $w(\lambda' ) =  \lambda'$ whenever 
$w(\lambda ) = \lambda$. Then $\pi (P,\delta ,\lambda' )$ is realized on the same 
vector space as $\pi (P,\delta ,\lambda )$, and it also decomposes as $V_1 \oplus V_2$.
\end{thm}

One deduces that Irr$(\mh H)$ is actually $\mf t / W$ with certain affine 
subspaces carrying a finite multiplicity. Since an affine space modulo a finite group is 
contractible, the cohomology of Irr$(\mh H)$ is rather easy, it is just a matter of 
counting affine subspaces modulo $W$, with multiplicities. All this leads to

\begin{thm}\label{thm:0.4}
Let $\mr{Irr}_0 (\mh H)$ be the collection of irreducible tempered 
$\mh H$-representations with real central character. Via the inclusion 
$\mh C [W] \to \mh H$ the set $\mr{Irr}_0 (\mh H)$ becomes a $\mh Q$-basis of
the representation ring $R(W) \otimes_{\mh Z} \mh Q$.
\end{thm}

Looking further ahead, via Lusztig's reduction theorems \cite{Lus-Gr} Theorem \ref{thm:0.4}
could have important consequences for the representation theory of affine Hecke
algebras. However, fur this purpose graded Hecke algebras are not quite sufficient,
one has to include automorphisms of the Dynkin diagram of $R$ in the picture. This 
prompted the author to generalize all the relevant representation theory to crossed
products of graded Hecke algebras with groups of diagram automorphisms.

Theorem \ref{thm:0.4} and part of the above are worked out in the sequel to
this paper \cite{Sol2}. Basically the author divided the material between these two papers
such that all the required representation theory is in this one, while the homological
algebra and cohomological computations are in \cite{Sol2}.
\vspace{2mm}

Let us briefly discuss the contents and the credits of the separate sections. In his attempts 
to provide a streamlined introduction to the representation theory of graded Hecke algebras 
the author acknowledges that this has already been done in \cite{KrRa} and \cite{Slo}. 
Nevertheless he found it useful to bring these and some other results together with a 
different emphasis, and meanwhile to fill in some gaps in the existing literature.

Section 1 contains the basic definitions of graded Hecke algebras. Section 2 describes 
parabolic induction and the Langlands classification, which can also be found in \cite{Eve,KrRa}. 
The normalized intertwining operators introduced in Section 3 appear to be new, 
although they look much like those in \cite{BaMo2} and those in \cite{Opd-Sp}. 
The results on affine Hecke algebras that we collect in Section 4 stem mostly from 
Opdam \cite{Opd-Sp}. In Section 5 we describe the link between affine and graded Hecke 
algebras in detail, building upon the work of Lusztig \cite{Lus-Gr}. In particular we determine
the global dimension of a graded Hecke algebra. In Section 6 we prove that tempered 
parabolically induced representations are unitary. This relies on work of Barbasch--Moy 
\cite{BaMo1,BaMo2} and Heckman--Opdam \cite{HeOp}. As already discussed above, we 
establish Theorems \ref{thm:0.1} and \ref{thm:0.2} in Section 7.
The generalization of the previous material to crossed products of graded Hecke algebras
with finite groups of diagram automorphisms is carried out in Section 8. Its representation
theoretic foundation is formed by Clifford theory, which we recall in the Appendix. 
Section 9 contains a more general version of Theorem \ref{thm:0.3}. Finally, in 
Section 10 we define a stratification on the spectrum of a graded Hecke algebra, and
we describe the strata as topological spaces.
\\[2mm]

\section{Graded Hecke algebras}
\label{sec:1}

For the construction of graded Hecke algebras we will use the following objects:
\begin{itemize}
\item a finite dimensional real inner product space $\mf a$,
\item the linear dual $\mf a^*$ of $\mf a$,
\item a crystallographic, reduced root system $R$ in $\mf a^*$,
\item the dual root system $R^\vee$ in $\mf a$,
\item a basis $\Pi$ of $R$.
\end{itemize}
We call 
\begin{equation}
\tilde{\mc R} = (\mf a^* ,R, \mf a, R^\vee, \Pi )
\end{equation}
a degenerate root datum. We do not assume that $\Pi$ spans
$\mf a^*$, in fact $R$ is even allowed to be empty.
Our degenerate root datum gives rise to
\begin{itemize}
\item the complexifications $\mf t$ and $\mf t^*$ of
$\mf a$ and $\mf a^*$,
\item the symmetric algebra $S (\mf t^*)$ of $\mf t^*$,
\item the Weyl group $W$ of $R$,
\item the set $S = \{ s_\alpha : \alpha \in \Pi \}$ of
simple reflections in $W$,
\item the complex group algebra $\mh C [W]$.
\end{itemize}
Choose formal parameters $\mb k_\alpha$ for $\alpha \in \Pi$,
with the property that $\mb k_\alpha = \mb k_\beta$ if
$\alpha$ and $\beta$ are conjugate under $W$. The graded
Hecke algebra $\tilde{\mh H} (\tilde{\mc R})$
corresponding to $\tilde{\mc R} (\tilde{\mc R})$ is defined as follows. As a
complex vector space
\[
\tilde{\mh H} = \mh C [W] \otimes S (\mf t^* ) \otimes 
\mh C [\{ \mb k_\alpha : \alpha \in \Pi \} ] .
\]
The multiplication in $\tilde{\mh H} (\tilde{\mc R})$ is determined by the following rules:
\begin{itemize}
\item $\mh C[W] \,, S (\mf t^* )$ and $\mh C[\{ \mb k_\alpha : \alpha \in \Pi \} ]$ 
are canonically embedded as subalgebras,
\item the $\mb k_\alpha$ are central in $\tilde{\mh H} (\tilde{\mc R})$,
\item for $x \in \mf t^*$ and $s_\alpha \in S$ we have the cross relation
\begin{equation}\label{eq:1.1}
x \cdot s_\alpha - s_\alpha \cdot s_\alpha (x) = 
\mb k_\alpha \inp{x}{\alpha^\vee} \,.
\end{equation}
\end{itemize}
We define a grading on $\tilde{\mh H} (\tilde{\mc R})$ by requiring that $\mf t^*$
and the $\mb k_\alpha$ are in degree one, while $W$ has degree zero.

In fact we will only study specializations of this algebra.
Pick complex numbers $k_\alpha \in \mh C$ for $\alpha \in \Pi$, such 
that $k_\alpha = k_\beta$ if $\alpha$ and $\beta$ are conjugate under
$W$. Let $\mh C_k$ be the onedimensional $\mh C[\{ \mb k_\alpha : \alpha \in \Pi \} ]$-module 
on which $\mb k_\alpha$ acts as multiplication by $k_\alpha$. We define 
\begin{equation}
\mh H = \mh H (\tilde{\mc R},k) = \tilde{\mh H} (\tilde{\mc R})
\otimes_{\mh C[\{ \mb k_\alpha : \alpha \in \Pi \} ]} \mh C_k .
\end{equation}
With some abuse of terminology $\mh H (\tilde{\mc R},k)$ is also
called a graded Hecke algebra. Notice that as a vector space
$\mh H (\tilde{\mc R},k)$ equals $\mh C[W] \otimes S (\mf t^*)$, and
that the cross relation \eqref{eq:1.1} now holds with $\mb k_\alpha$
replaced by $k_\alpha$:
\begin{equation}\label{eq:1.2}
x \cdot s_\alpha - s_\alpha \cdot s_\alpha (x) = k_\alpha \inp{x}{\alpha^\vee} .
\end{equation}
Since $S (\mf t^*)$ is Noetherian and $W$ is finite, $\mh H$ is Noetherian as well.
We define a grading on $\mh H$ by deg$(w) = 0 \; \forall w \in W$ and
deg$(x) = 1 \; \forall x \in \mf t^*$. However, the algebra $\mh H'$ is in general not 
graded, only filtered. That is, the product $h_1 h_2$ of two homogeneous elements 
$h_1 ,h_2 \in \mh H'$ need not be homogeneous, but all its homogeneous components
have degree at most $\text{deg}(h_1 ) + \text{deg}(h_2 )$.
Let us mention some special cases in which $\mh H (\tilde{\mc R},k)$ is graded:
\begin{itemize}
\item if $R = \es$ then $\mh H = \tilde{\mh H} (\tilde{\mc R}) = S (\mf t^* )$,
\item if $k_\alpha = 0 \; \forall \alpha \in \Pi$ then $\mh H (\tilde{\mc R},k)$ is 
the usual crossed product $W \ltimes S (\mf t^* )$, with the cross relations
\begin{equation}
w \cdot x = w(x) \cdot w \qquad w \in W , x \in \mf t^* .
\end{equation}
\end{itemize}
More generally, for any $s_\alpha \in S$ and $p \in S (\mf t^* )$ we have
\[
p \cdot s_\alpha - s_\alpha \cdot s_\alpha (p) = k_\alpha \frac{p - s_\alpha (p)}{\alpha} .
\]
Multiplication with any $z \in \mh C^\times$ defines a bijection $m_z : \mf t^* \to \mf t^*$,
which clearly extends to an algebra automorphism of $S(\mf t^* )$. From the cross relation
\eqref{eq:1.2} we see that it extends even further, to an algebra isomorphism
\begin{equation}\label{eq:1.3}
m_z : \mh H (\tilde{\mc R},zk) \to \mh H (\tilde{\mc R}, k)
\end{equation}
which is the identity on $\mh C[W]$. In particular, if all $\alpha \in R$ are conjugate 
under $W$, then there are essentially only two graded Hecke algebras attached to 
$\tilde{\mc R}$: one with $k=0$ and one with $k \neq 0$.
\\[2mm]

\section{Parabolic induction}
\label{sec:2}

A first tool to study $\mh H$-modules is restriction to the commutative
subalgebra $S (\mf t^* ) \subset \mh H$. Let $(\pi ,V)$ be an $\mh H$-module and pick
$\lambda \in \mf t$. The $\lambda$-weight space of $V$ is
\[
V_\lambda = \{ v \in V : \pi (x) v = \inp{x}{\lambda} v \; \forall x \in \mf t^* \} \,,
\] 
and the generalized $\lambda$-weight space is
\[
V_\lambda^{gen} = \{ v \in V : \exists n \in \mh N :
(\pi (x) - \inp{x}{\lambda} )^n v = 0 \, \forall x \in \mf t^* \} \,.
\]
We call $\lambda$ a $S (\mf t^* )$-weight of $V$ if $V_\lambda^{gen} \neq 0$,
or equivalently if $V_\lambda \neq 0$. If $V$ has finite dimension then it decomposes as
\begin{equation}
V = {\ts \bigoplus_{\lambda \in \mf t}} V_\lambda^{gen} \,.
\end{equation}
According to \cite[Proposition 4.5]{Lus-Gr} the center of $\mh  H$ is
\begin{equation}
Z (\mh H ) = S (\mf t^* )^W .
\end{equation}
In particular $\mh H$ is of finite rank as a module over its center, so all its irreducible
modules have finite dimension. Moreover the central character of an irreducible 
$\mh H$-module can be regarded as an element of $\mf t / W$.

Let $P \subset \Pi$ be a set of simple roots. They form a basis of a root subsystem
$R_P \subset R$ with Weyl group $W_P \subset W$. Let $\mf a_P \subset \mf a$ and 
$\mf a_P^* \subset \mf a^*$ be the real spans of respectively $R_P^\vee$ and $R_P$.
We denote the complexifications of these vector spaces by $\mf t_P$ and $\mf t_P^*$, 
and we write 
\[
\begin{array}{lllll}
\mf t^P & = & (\mf t_P^* )^\perp & = &
\{ \lambda \in \mf t : \inp{x}{\lambda} = 0 \: \forall x \in \mf t_P^* \} \,, \\
\mf t^{P*} & = & (\mf t_P )^\perp & = & 
\{ x \in \mf t^* : \inp{x}{\lambda} = 0 \: \forall \lambda \in \mf t_P \} \,.
\end{array}
\]
We define the degenerate root data 
\begin{align}
& \tilde{\mc R}_P = (\mf a_P^* ,R_P ,\mf a_P ,R_P^\vee ,P) \,, \\
& \tilde{\mc R}^P = (\mf a^* ,R_P ,\mf a ,R_P^\vee ,P) \,,
\end{align}
and the graded Hecke algebras
\begin{align}
\mh H_P = \mh H (\tilde{\mc R}_P ,k) \,, \\
\mh H^P = \mh H (\tilde{\mc R}^P ,k) \,.
\end{align}
Notice that the latter decomposes as a tensor product of algebras: 
\begin{equation}
\mh H^P = S (\mf t^{P*} ) \otimes \mh H_P .
\end{equation}
In particular every irreducible $\mh H^P$-module is of the form $\mh C_\lambda \otimes V$, 
where $\lambda \in \mf t^P$ and $V$ is an irreducible $\mh H_P$-module.
In general, for any $\mh H_P$-module $(\rho, V_\rho )$ and $\lambda \in \mf t^P$
we denote the action of $\mh H^P$ on $\mh C_\lambda \otimes V_\rho$ by 
$\rho_\lambda$. We define the parabolically induced module
\begin{equation}
\pi (P,\rho ,\lambda ) = \mr{Ind}_{\mh H^P}^{\mh H} (\mh C_\lambda \otimes V_\rho )
= \mr{Ind}_{\mh H^P}^{\mh H} (\rho_\lambda ) = \mh H \otimes_{\mh H^P} V_{\rho_\lambda} \,.
\end{equation}
We remark that these are also known as ``standard modules'' \cite{BaMo1,Kri}.
Of particular interest is the case $P = \emptyset$. Then $\mh H^P = S(\mf t^* ) \,, 
\mh H_P = \mh C$ and we denote the unique irreducible representation of $\mh H_{\emptyset}$ 
by $\delta_0$. The parabolically induced $\mh H$-modules $\pi (\emptyset ,\delta_0 ,\lambda )$ 
form the principal series, which has been studied a lot. For example, it is easily shown 
that every irreducible $\mh H$-representation is a quotient of some principal series 
representation \cite[Proposition 4.2]{Kri}. Let
\begin{equation}
W^P := \{ w \in W : \ell (w s_\alpha ) > \ell (w) \; \forall \alpha \in P \}
\end{equation}
be the set of minimal length representatives of $W / W_P$.

\begin{lem}\label{lem:2.2} \textup{\cite[Theorem 6.4]{BaMo2}}\\
The weights of $\pi (P,\rho ,\lambda )$ are precisely the elements $w (\lambda + \mu )$, 
where $\mu$ is a $S(\mf t_P^* )$-weight of $\rho$ and $w \in W^P$.
\end{lem}

Since the complex vector space $\mf t$ has a distinguished real form $\mf a$, we can
decompose any $\lambda \in \mf t$ unambiguously as
\begin{equation}
\lambda = \Re (\lambda ) + i \Im (\lambda ) \text{ with } 
\Re (\lambda ), \Im (\lambda ) \in \mf a \,.
\end{equation}
We define the positive cones 
\begin{equation}
\begin{array}{lll}
\mf a^{*+} & = & \{ x \in \mf a^* : 
  \inp{x}{\alpha^\vee} \geq 0 \: \forall \alpha \in \Pi \} \,, \\
\mf a_P^+ & = &  \{ \mu \in \mf a_P : 
  \inp{\alpha}{\mu} \geq 0 \: \forall \alpha \in P \} \,, \\
\mf a^{P+} & = & \{ \mu \in \mf a^P : 
  \inp{\alpha}{\mu} \geq 0 \: \forall \alpha \in \Pi \setminus P \} \,, \\
\mf a^{P++} & = & \{ \mu \in \mf a^P : 
  \inp{\alpha}{\mu} > 0 \; \forall \alpha \in \Pi \setminus P \} \,.
\end{array}
\end{equation}
The antidual of $\mf a^{*+}$ is
\begin{equation} 
\mf a^- = \{ \lambda \in \mf a : \inp{x}{\lambda} \leq 0 \: \forall x \in \mf a^{*+} \} = 
\big\{ {\ts \sum_{\alpha \in \Pi}} \lambda_\alpha \alpha^\vee : \lambda_\alpha \leq 0 \big\} \,.
\end{equation}
The interior $\mf a^{--}$ of $\mf a^-$ equals
\[
\big\{ {\ts \sum_{\alpha \in \Pi}} \lambda_\alpha \alpha^\vee : \lambda_\alpha < 0 \big\}
\]
if $\Pi$ spans $\mf a^*$, and is empty otherwise.
A finite dimensional $\mh H$-module $V$ is called tempered if $\Re (\lambda ) \in \mf a^-$,
for all weights $\lambda$. More restrictively we say that $V$ belongs to the
discrete series if it is irreducible and $\Re (\lambda ) \in \mf a^{--}$, again for all
weights $\lambda$.

\begin{lem}\label{lem:2.3}
Let $\rho$ be a finite dimensional $\mh H_P$-module and let $\lambda \in \mf t^P$.
The $\mh H$-representation $\pi (P,\rho,\lambda )$ is tempered if and only if $\rho$
is tempered and $\lambda \in i \mf a^P$.
\end{lem}
\emph{Proof.}
If $\rho$ is tempered and $\lambda \in i \mf a^P$, then $\rho_\lambda$ is  a tempered
$\mh H^P$-representation. According to \cite[Corollary 6.5]{BaMo1} $\pi (P,\rho,\lambda )$
is a tempered $\mh H$-representation.

Conversely, suppose that $\pi (P,\rho ,\lambda )$ is tempered, and let $\mu$ be any 
$S(\mf t^* )$-weight of $\rho_\lambda$. By Lemma \ref{lem:2.2} $w(\mu )$ is  a weight of 
$\pi (P,\rho ,\lambda )$, for every $w \in W^P$. As is well known \cite[Section 1.10]{Hum}
\[
W^P = \{ w \in W : \ell (w s_\alpha ) > \ell (w) \; \forall \alpha \in P \} =
\{ w\in W : w(P) \subset R^+ \} \,.
\]
We claim that $W^P$ can be characterized in the following less obvious way:
\begin{equation}\label{eq:2.1}
{\ts \bigcup_{w \in W^P}} \, w^{-1} (\mf a^{*+} ) = 
\{ x \in \mf a^* : \inp{x}{\alpha^\vee} \geq 0 \; \forall \alpha \in P \} \,.
\end{equation}
Since $P \subset w^{-1}(R^+ ) \: \forall w \in W^P$, the inclusion $\subset$ holds. 
The right hand side is the positive chamber for the root system $R_P$ in $\mf a^*$, 
so it is a fundamental domain for action of $W_P$ on $\mf a^*$. Because
\[
\bigcup_{w \in W^P} W^P w^{-1} = \bigcup_{w \in W^P} w W_P = W \,, 
\]
the left hand side also intersects every $W_P$-orbit in $\mf a^*$. Thus both sides of
\eqref{eq:2.1} are indeed equal. Taking antiduals we find
\begin{multline}
\{ \nu \in \mf a : \inp{x}{w (\nu )} \leq 0 \; \forall x \in \mf a^{*+} , w \in W^P \} \\
= \{ \nu \in \mf a : \inp{x}{\nu} \leq 0 \; 
\forall x \in {\ts \bigcup_{w \in W^P}} w^{-1} (\mf a^{*+} ) \} = \mf a_P^- \,.
\end{multline}
Combining this with the definition of temperedness, we deduce that $\Re (\mu ) \in 
\mf a_P^-$. Since $\lambda \in \mf t^P$ and $\mu - \lambda \in \mf t_P$, we conclude
that $\Re (\lambda ) = 0$ and $\lambda \in i \mf a^P$. Furthermore we see now that every
weight $\mu - \lambda$ of $\rho$ has real part in $\mf a_P^-$, so $\rho$ is tempered.
$\qquad \Box$ \\[2mm]

The Langlands classification explains how to reduce the classification of irreducible
$\mh H$-modules to that of irreducible tempered modules. We say that a triple
$(P,\rho ,\nu )$ is a Langlands datum if
\begin{itemize}
\item $P \subset \Pi$,
\item $(\rho ,V_\rho )$ is an irreducible tempered $\mh H_P$-module,
\item $\nu \in \mf t^P$,
\item $\Re (\nu ) \in \mf a^{P++}$.
\end{itemize}

\begin{thm}\label{thm:2.1}
\textup{(Langlands classification)}
\begin{description}
\item[a)] For every Langlands datum $(P,\rho,\nu )$ the $\mh H$-module 
$\pi (P,\rho ,\nu ) = \mr{Ind}_{\mh H^P}^{\mh H} (\mh C_\nu \otimes V_\rho )$
has exactly one irreducible quotient, which we call $L(P,\rho ,\nu )$.
\item[b)] If $(Q,\sigma ,\mu )$ is another Langlands datum and
$L(Q,\sigma ,\mu )$ is equivalent to $L(P,\rho ,\nu )$, then $Q = P \,,
\mu = \nu$ and $\sigma \cong \rho$.
\item[c)] For every irreducible $\mh H$-module $V$ there exists a Langlands datum
such that\\ $V \cong L (P,\rho,\nu)$.
\end{description}
\end{thm}
\emph{Proof.}
See \cite{Eve} or \cite[Theorem 2.4]{KrRa}. $\qquad \Box$
\\[2mm]
We would like to know a little more about the relation between the Langlands quotient 
and the other constituents of $\pi (P,\rho ,\nu )$. The proof shows that $L (P,\rho ,\nu )$
has a highest weight and is cyclic for $\pi (P,\rho ,\nu )$, in a suitable sense. These
properties are essential in the following result.

Denote the central character of any irreducible $\mh H_P$-module $\delta$ by
\begin{equation}
cc_P (\delta ) \in \mf t_P / W_P \,,
\end{equation}
and identify it with the corresponding $W_P$-orbit in $\mf t_P$. Since $W_P$ acts 
orthogonally on $\mf a_P$, the number $\| \Re (\lambda ) \|$ is the same for all 
$\lambda \in cc_P (\delta )$, and hence may be written as $\| \Re (cc_P (\delta )) \|$.

\begin{prop}\label{prop:2.4}
Let $(P,\rho ,\nu )$ and $(Q,\sigma ,\mu )$ be different Langlands data,
and let $(\rho' ,V')$ be another irreducible tempered $\mh H_P$-module.
\begin{description}
\item[a)] The functor $\mr{Ind}_{\mh H^P}^{\mh H}$ induces an isomorphism
\[
\mr{Hom}_{\mh H_P} (\rho ,\rho' ) = \mr{Hom}_{\mh H^P} (\rho_\nu ,\rho'_\nu )
\cong \mr{Hom}_{\mh H} (\pi (P,\rho ,\nu ) ,\pi (P,\rho' ,\nu )) \,.
\]
In particular these spaces are either 0 or onedimensional.
\item[b)] Suppose that $L(Q,\sigma ,\mu )$ is a constituent of $\pi (P,\rho ,\nu )$. Then 
\[
P \subset Q \quad and \quad \| \Re (cc_P (\rho )) \| < \| \Re (cc_Q (\sigma )) \| \,.
\]
\end{description}
\end{prop}
\emph{Proof.}
a) Since $S (\mf t^{P*}) \subset \mh H^P$ acts on both $\rho_\nu$ and $\rho'_\nu$
by the character $\nu$, we have
\[ 
\mr{Hom}_{\mh H_P} (\rho ,\rho' ) = \mr{Hom}_{\mh H^P} (\rho_\nu ,\rho'_\nu ) \,.
\]
For $\alpha \in \Pi$ we define $\delta_\alpha \in \mf a_\Pi$ by
\[
\inp{\beta}{\delta_\alpha} = \left\{ \begin{array}{lll}
1 & \mr{if} & \alpha = \beta \\
0 & \mr{if} & \alpha \neq \beta \in \Pi \,.
\end{array} \right.
\]
According to Langlands \cite[Lemma 4.4]{Lan}, for every $\lambda \in \mf a$
there is a unique subset $F (\lambda ) \subset \Pi$ such that $\lambda$ can
be written as
\begin{equation}\label{eq:2.3}
\lambda = \lambda^\Pi + \sum_{\alpha \in \Pi \setminus F(\lambda )} d_\alpha \delta_\alpha
+ \sum_{\alpha \in F (\lambda )} c_\alpha \alpha^\vee \quad
\mr{with} \; \lambda^\Pi \in a^\Pi , d_\alpha > 0, c_\alpha \leq 0 \,.
\end{equation}
We put $\lambda_0 = \sum_{\alpha \in \Pi \setminus F(\lambda )} d_\alpha \delta_\alpha 
\in \mf a^+$. For any weight $\lambda$ of $\rho_\nu$ we have $\Re (\lambda -\nu ) 
\in \mf a_P^-$ and $(\Re \lambda )_0 = \Re \nu |_{\mf t^*_\Pi }$.

Let $\lambda'$ be a weight of $\rho'_\nu$. According to \cite[(2.13)]{KrRa}
\begin{equation}\label{eq:2.2}
\Re (w \lambda' )_0 < \Re (\lambda' )_0 = \Re \nu |_{\mf t^*_\Pi } \qquad
\forall w \in W^P \setminus \{ 1 \} \,,
\end{equation}
with respect to the ordering that $\Pi$ induces on $\mf a_\Pi^*$.
Hence for $w \in W^P , w(\lambda' )$ can only equal $\lambda$ if $w = 1$.
Let $v_\lambda \in \mh C_\nu \otimes V_\rho$ be a nonzero weight vector. Since 
$\rho_\lambda$ is an irreducible $\mh H^P$-representation, $1 \otimes v_\lambda \in 
\mh H \otimes_{\mh H^P} \mh C_\nu \otimes V_\rho$ is cyclic for $\pi (P,\rho ,\nu )$, and
therefore the map
\[
\mr{Hom}_{\mh H} (\pi (P,\rho ,\nu ) ,\pi (P,\rho' ,\nu )) \to \pi (P,\rho' ,\nu ) : 
f \mapsto f (1 \otimes v_\lambda )
\]
is injective. By \eqref{eq:2.2} the $\lambda$-weight space of $\pi (P,\rho' ,\nu )$
is contained in $1 \otimes \mh C_\nu \otimes V'$. (This weight space might be zero. 
See also the more general calculations on page \pageref{eq:9.2}.)
So $f (1 \otimes v_\lambda ) \in 1 \otimes \mh C_\nu \otimes V'$ and in fact
\[
f (1 \otimes \mh C_\nu \otimes V_\rho ) \subset 1 \otimes \mh C_\nu \otimes V' \,.
\]
Thus any $f \in \mr{Hom}_{\mh H} (\pi (P,\rho ,\nu ) ,\pi (P,\rho' ,\nu ))$ lies in
$\mr{Ind}_{\mh H^P}^{\mh H} \big( \mr{Hom}_{\mh H^P} (\rho_\nu ,\rho'_\nu ) \big)$.\\
b) Since $S \big( \mf t^{\Pi *} \big)$ acts on $\pi (P,\rho ,\nu )$ by $\nu$ and on 
$\pi (Q,\sigma ,\mu )$ by $\mu$, we have $\nu |_{\mf t^{\Pi *}} = 
\mu |_{\mf t^{\Pi *}}$. Therefore $S \big( \mf t^{\Pi *} \big)$ 
presents no problems, and we may just as well assume that $\nu ,\mu \in \mf t_\Pi^*$.

By construction \cite[p. 39]{KrRa} $L(P,\rho ,\nu )$ is the unique irreducible 
subquotient of $\pi (P.\rho ,\nu )$ which has a $S(\mf t^* )$-weight $\lambda$ with
$(\Re \lambda )_0 = \Re \nu$. Moreover of all weights $\lambda'$ of proper submodules of
$\pi (P,\rho ,\nu )$ satisfy $(\Re \lambda' )_0 < \Re \nu$, with the notation 
of \eqref{eq:2.3}. In particular, for the
subquotient $L(Q, \sigma ,\mu )$ of $\pi (P,\rho ,\nu )$ we find that $\Re \mu < \Re \nu$.
Since $\Re \mu \in \mf a^{Q++}$ and $\Re \nu \in \mf a^{P++}$, this implies $P \subset Q$ and
$\| \Re \mu \| < \| \Re \nu \|$.

According to Lemma \ref{lem:2.2} all constituents of $\pi (P,\rho ,\nu )$ have central 
character $W (cc_P (\rho ) + \nu ) \in \mf t / W$. The same goes for $(Q,\sigma ,\mu )$, so
\[
W (cc_P (\rho ) + \nu ) = W (cc_Q (\sigma ) + \mu ) \,.
\]
By definition $\nu \perp \mf t_P$ and $\mu \perp \mf t_Q$, so
\[
\| \Re (cc_P (\rho )) \|^2 + \| \Re \nu \|^2 = \| \Re (cc_P (\rho ) + \nu) \|^2 = 
\| \Re (cc_Q (\sigma ) + \mu ) \|^2 = \| \Re (cc_Q (\sigma )) \|^2 + \| \mu \|^2 .
\]
Finally we use that $\| \Re \mu \|^2 < \| \Re \nu \|^2 . \qquad \Box$
\vspace{4mm}

\section{Intertwining operators}
\label{sec:4} 

We will construct rational intertwiners between parabolically induced representations.
Our main ingredients are the explicit calculations of Lusztig \cite{Lus-Gr} and 
Opdam's method for constructing the corresponding intertwiners in the context of
affine Hecke algebras. 

Let $\mh C (\mf t /W) = \mh C (\mf t )^W$ be the quotient field of 
\[
\mh C [\mf t /W] = \mh C [\mf t]^W = S (\mf t^* )^W = Z (\mh H ) \,,
\] 
and write
\[
\prefix{_\mc F}{}{\mh H} = \mh C (\mf t )^W \otimes_{Z (\mh H)} \mh H = 
\mh C (\mf t ) \otimes_{S (\mf t^*)} \mh H \,.
\]
For $\alpha \in \Pi$ we define 
\[
\tilde \tau_{s_\alpha} := (s_\alpha \alpha - k_\alpha )(\alpha + k_\alpha )^{-1}
\in \prefix{_\mc F}{}{\mh H} \,.
\]
\begin{prop}\label{prop:4.1}
\begin{description} The elements $\tilde \tau_{s_\alpha}$ have the following properties:
\item[a)] The map $s_\alpha \mapsto \tilde \tau_{s_\alpha}$ extends to a group
homomorphism $W \to ( \prefix{_\mc F}{}{\mh H} )^\times$.
\item[b)] For $f \in \mh C (\mf t )$ and $w \in W$ we have 
$\tilde \tau_w f \tilde \tau_w^{-1} = w (f)$.
\item[c)] The map
\[
\begin{array}{ccc}
\mh C (\mf t ) \rtimes W & \to & \prefix{_\mc F}{}{\mh H} \\
\sum_{w \in W} f_w w & \mapsto & \sum_{w \in W} f_w \tilde \tau_w
\end{array}
\]
is an algebra isomorphism.
\end{description}
\end{prop}
\emph{Proof.}
See \cite[Section 5]{Lus-Gr}. $\qquad \Box$
\\[3mm]

As this proposition already indicates, conjugation by $\tilde \tau_w$ will
be a crucial operation. For $P,Q \subset \Pi$ we define
\[
W(P,Q) = \{w \in W : w(P) = Q \} \,.
\]
\begin{lem}\label{lem:4.2}
\begin{description}
\item[a)] Let $P,Q \subset \Pi , u \in W_P$ and $w \in W(P,Q)$.
Then $\tilde \tau_w u \tilde \tau_w^{-1} = w u w^{-1}$.
\item[b)] There are algebra isomorphisms
\begin{align*}
& \psi_w : \mh H_P \to \mh H_Q \\
& \psi_w : \mh H^P \to \mh H^Q \\
& \psi_w (x u) = \tilde \tau_w x u \tilde \tau_w^{-1} = w(x) \, w u w^{-1} 
\qquad x \in \mf t^* , u \in W_P \,.
\end{align*}
\end{description}
\end{lem}
\emph{Proof.}
First we notice that b) is an immediate consequence of a) and Proposition 
\ref{prop:4.1}.b. It suffices to show a) for $u = s_\alpha$ with $\alpha \in P$.
Instead of a direct calculation, our strategy is to show that the algebra 
homomorphism $f \to w(f) = \tilde \tau_w f \tilde \tau_w^{-1}$ has only one
reasonable extension to $W_P$. Pick $x \in \mf t^*$ and write 
$\beta = w(\alpha ) \in \Pi$. By Proposition \ref{prop:4.1}.b 
\begin{equation}
\begin{array}{lll}
k_\alpha \inp{x}{\alpha^\vee} & = & 
 \tilde \tau_w k_\alpha \inp{x}{\alpha^\vee} \tilde \tau_w^{-1} \\
 & = & \tilde \tau_w (x s_\alpha - s_\alpha s_\alpha (x) ) \tilde \tau_w^{-1} \\
 & = & w(x) \tilde \tau_w s_\alpha \tilde \tau_w^{-1} - 
 \tilde \tau_w s_\alpha \tilde \tau_w^{-1} s_\beta (x) \,. 
\end{array}
\end{equation}
On the other hand
\[
w(x) s_\beta - s_\beta (s_\beta w)(x) = k_\beta \inp{w(x)}{\beta^\vee} = 
k_\alpha \inp{x}{\alpha^\vee} \,.
\]
So for every $y = w(x) \in \mf t^*$ we have
\begin{equation}\label{eq:4.1}
y (\tilde \tau_w s_\alpha \tilde \tau_w^{-1} - s_\beta) =
(\tilde \tau_w s_\alpha \tilde \tau_w^{-1} - s_\beta ) s_\beta (y) \,.
\end{equation}
Using Proposition \ref{prop:4.1}.c we can write
\[
\tilde \tau_w s_\alpha \tilde \tau_w^{-1} - s_\beta =
{\ts \sum_{v \in W}} \tilde \tau_v f_v \qquad f_v \in \mh C (\mf t )
\]
in a unique way. Comparing \eqref{eq:4.1} with the multiplication in
$\mh C (\mf t ) \rtimes W$ we find that $f_v = 0$ for $v \neq s_\beta$, so
\[
\tilde \tau_w s_\alpha \tilde \tau_w^{-1} = s_\beta + \tilde \tau_{s_\beta} f_{s_\beta}
= \tilde \tau_{s_\beta} \left( f_{s_\beta} + (k_\beta + \beta) \beta^{-1} \right) + 
k_\beta \beta^{-1} = \tilde \tau_{s_\beta} f + k_\beta \beta^{-1} \,,
\]
with $f = f_{s_\beta} + (k_\beta + \beta ) \beta^{-1} \in \mh C (\mf t )$. But
\[
\begin{array}{lll}
1 = s_\alpha^2 & = & (\tilde \tau_w s_\alpha \tilde \tau_w^{-1} )^2 \\
 & = & ( \tilde \tau_{s_\beta} f + k_\beta \beta^{-1} )^2 \\
 & = & \tilde \tau^2_{s_\beta} s_\beta (f) f + \tilde \tau_{s_\beta} f k_\beta \beta^{-1} +
 k_\beta \tilde \tau_{s_\beta} s_\beta (\beta^{-1}) f + k_\beta^2 \beta^{-2} \\
 & = & s_\beta (f) f + k_\beta \tilde \tau_{s_\beta} \beta^{-1} f -
 k_\beta \tilde \tau_{s_\beta} \beta^{-1} f + k_\beta^2 \beta^{-2} \\
 & = & s_\beta (f) f + k_\beta^2 \beta^{-2}
\end{array}
\]
Writing $f = f_1 / f_2$ with $f_i \in \mh C [\mf t]$ we find
\begin{align*}
& s_\beta (f) f = 1 - k_\beta^2 \beta^{-2} = \frac{\beta^2 - k_\beta^2}{\beta^2} = 
s_\beta \left( \frac{k_\beta + \beta}{\beta} \right) \frac{k_\beta + \beta}{\beta} \,, \\
& s_\beta (f_1 \beta ) f_1 \beta = s_\beta ( f_2 (k_\beta + \beta )) f_2 (k_\beta + \beta ) \,,
\end{align*}
which is only possible if $f_1 \beta = \pm f_2 (k_\beta + \beta )$. 
Equivalently either $f_{s_\beta} = 0$ or \\
$f_{s_\beta} = -2(k_\beta + \beta)\beta^{-1}$, and either
\begin{align}
& \label{eq:4.2} \tilde \tau_w s_\alpha \tilde \tau_w^{-1} = s_\beta \qquad \text{or} \\
& \nonumber \tilde \tau_w s_\alpha \tilde \tau_w^{-1} = 
 s_\beta - 2 \tilde \tau_{s_\beta} (k_\beta + \beta)\beta^{-1} = s_\beta - 2(s_\beta + 1) +
 2 (k_\beta + \beta)\beta^{-1} = 2 k_\beta \beta^{-1} - s_\beta \,.
\end{align}
However, all the above expressions are rational in the parameters $k$, and for $k = 0$ we
clearly have $\tilde \tau_w s_\alpha \tilde \tau_w^{-1} = w s_\alpha w^{-1} = s_\beta$.
Hence the second alternative of \eqref{eq:4.2} cannot hold for any $k. \qquad \Box$
\\[2mm]

As above, let $w \in W(P,Q)$ and take $\lambda \in \mf t^P$.
Let $(\rho ,V_\rho )$ be any finite dimensional $\mh H_P$-module, and let 
$(\sigma ,V_\sigma )$ be a $\mh H_Q$-module which is equivalent to $\rho \circ \psi_w^{-1}$. 
Our goal is to construct an intertwiner between the $\mh H$-modules $\pi (P,\rho,\lambda)$ 
and $\pi (Q,\sigma,w(\lambda ))$. The (nonnormalized) intertwining operators from 
\cite[Section 1.6]{BaMo2} are insufficient for our purposes, as they do not match up with
the corresponding (normalized) intertwiners for affine Hecke algebras. This requires the
use of $\tilde{\tau}_{s_\alpha}$ and not just $s_\alpha \alpha - k_\alpha$. By assumption 
there exists a linear bijection $I_\rho^w : V_\rho \to V_\sigma$ such that
\begin{equation}\label{eq:4.3}
I_\rho^w (\rho_\lambda (h) v) = \sigma_{w (\lambda )} (\psi_w h) (I_\rho^w v) 
\qquad \forall h \in \mh H^P , v \in V_\rho \,.
\end{equation}
Consider the map
\begin{equation}\label{eq:4.5}
\begin{split}
& I_w : \prefix{_\mc F}{}{\mh H} \otimes_{\mh H^P} (\mh C_\lambda \otimes V_\rho ) \to
\prefix{_\mc F}{}{\mh H} \otimes_{\mh H^Q} (\mh C_{w (\lambda )} \otimes V_\sigma ) \\
& I_w (h \otimes_{\mh H^P} v) = h \tilde \tau_w^{-1} \otimes_{\mh H^Q} I_\rho^w (v)
\end{split}
\end{equation}
We check that it is well-defined:
\[
\begin{array}{lll}
I_w (h \otimes \rho_\lambda (h') v) & = & 
h \tilde \tau_w^{-1} \otimes I_\rho^w (\rho_\lambda (h') v) \\
 & = & h \tilde \tau_w^{-1} \otimes \sigma_{w (\lambda )} (\psi_w h') (I_\rho^w v) \\
 & = & h \tilde \tau_w^{-1} \psi_w (h') \otimes I_\rho^w (v) \\
 & = & h \tilde \tau_w^{-1} \tilde \tau_w h' \tilde \tau_w^{-1} \otimes I_\rho^w (v)
 \quad = \quad I_w (h h' \otimes v)
\end{array}
\]
Notice that due to some freedom in the construction, $I_u \circ I_w$ need not equal 
$I_{uw} \; (u,w \in W)$. Let $s_1 \cdots s_r$ be a reduced expression for 
$w^{-1} \in W$, with $s_i = s_{\alpha_i}$ simple reflections. 
\[
\begin{split}
& \tilde \tau_w^{-1} \;=\; \tilde \tau_{s_1} \cdots \tilde \tau_{s_r} \;=\;  
(s_1 \alpha_1 - k_1 )(\alpha_1 + k_1 )^{-1} \cdots 
 (s_r \alpha_r - k_r )(\alpha_r + k_r )^{-1} \;= \\
& (s_1 \alpha_1 - k_1 ) \cdots (s_r \alpha_r - k_r ) (s_r \cdots s_2 (\alpha_1 ) + k_1 )^{-1}
 \cdots (s_r (\alpha_{r-1}) + k_{r-1})^{-1} (\alpha_r + k_r )^{-1} \\
& =\; \prod_{i=1}^r (s_i \alpha_i - k_i ) 
 \prod_{\alpha \in R^+ : w^{-1}(\alpha )\in R^-} (\alpha + k_\alpha )^{-1}
\end{split}
\]
For any $S (\mf t_Q^* )$-weight $\mu$ of $\sigma$, the function 
$\prod_{\alpha \in R^+ : w^{-1}(\alpha )\in R^-} (\alpha + k_\alpha )^{-1}$
is regular on a nonempty Zariski-open subset of $\mu + \mf t^Q$, because
$w^{-1}(Q) = P \subset R^+$. Since $\sigma$ has only finitely many weights, this implies that
\[
\sigma_\nu \Big( {\ts \prod_{\alpha \in R^+ : w^{-1} (\alpha )\in R^-}} 
(\alpha + k_\alpha )^{-1} \Big)
\]
is invertible for all $\nu$ in a certain Zariski-open $U \subset \mf t^Q$. Hence $\sigma_\nu$
extends to a representation of a suitable algebra containing $\mh H^Q$ and $\tilde \tau_w^{-1}$.
Moreover the map 
\[
\mf t^Q \to V_\sigma : \nu \mapsto \sigma_\nu (\tilde \tau_w^{-1}) v
\]
is rational, with poles exactly in $\mf t^Q \setminus U$. So for $\nu \in U \,, I_w$ 
restricts to a map
\begin{equation}\label{eq:4.6}
\pi (w,P,\rho,\lambda ) : \mh H \otimes_{\mh H^P} (\mh C_\lambda \otimes V_\rho ) \to
\mh H \otimes_{\mh H^Q} (\mh C_{w (\lambda )} \otimes V_\sigma ) \,.
\end{equation}

\begin{prop}\label{prop:4.3}
The intertwining operator $\pi (w,P,\rho,\lambda )$ is rational as a function of 
$\lambda \in \mf t^P$. It is regular and invertible on a dense Zariski-open subset of $\mf t^P$.
\end{prop}
\emph{Proof.}
Everything except the invertibility was already discussed. Clearly
\[
I_w^{-1}(h \otimes v) = h \tilde \tau_w \otimes \big( I_\rho^w \big)^{-1} (v) = 
I_{w^{-1}} (h \otimes v) \,.
\]
By the same reasoning as above, the operator $I_{w^{-1}}$ is regular for $\lambda$ in a 
nonempty Zariski-open subset of $\mf t^P. \qquad \Box$
\\[2mm]

We remark that it is usually hard to determine $\pi (w,P,\delta ,\lambda )$ explicitly,
at least if dim $V_\rho > 1$.
Although in general $\pi (w,P,\rho,\lambda )$ cannot be extended continuously to all 
$\lambda \in \mf t^P$, we can nevertheless draw some conclusions that hold for all 
$\lambda \in \mf t^P$.

\begin{lem}\label{lem:4.4}
The $\mh H$-modules $\pi (P,\rho,\lambda )$ and $\pi (Q,\sigma,w(\lambda ))$ have the same
irreducible constituents, counted with multiplicity.
\end{lem}
\emph{Proof.}
Since $\mh H$ is of finite rank over its center, the Frobenius--Schur Theorem 
(cf. \cite[Theorem 27.8]{CuRe}) tells us that the characters of inequivalent irreducible
$\mh H$-modules are linearly independent functionals. Hence the lemma is equivalent to
\begin{equation}\label{eq:4.4}
\text{tr } \pi (P,\rho,\lambda )(h) - \text{tr } \pi (Q,\sigma,w(\lambda )) (h) = 0 
\qquad \forall h \in \mh H \,.
\end{equation}
By Proposition \ref{prop:4.3} we have $\pi (P,\rho,\lambda )(h) \cong 
\pi (Q,\sigma,w(\lambda ))$ for $\lambda$ in a Zariski-dense subset of $\mf t^P$. Hence the
left hand side of \eqref{eq:4.4} is 0 on a dense subset of $\mf t^P$. Finally we note that 
for fixed $h \in \mh H$, it is a polynomial function of $\lambda \in \mf t^P. \qquad \Box$
\\[2mm]

\section{Affine Hecke algebras}
\label{sec:3}

We will introduce the most important objects in the theory of affine Hecke algebras.
As far as possible we will use the notations from Section \ref{sec:1}. Most of the
things that we will claim can be found in \cite{Opd-Sp}.

Let $Y \subset \mf a$ be a lattice, and let $X = \mr{Hom}_{\mh Z} (Y,\mh Z ) \subset 
\mf a^*$ be its dual lattice. We assume that $R \subset X$ and $R^\vee \subset Y$.
Thus we have a based root datum
\[
\mc R = (X,R,Y,R^\vee ,\Pi ) \,.
\]
We are interested in the affine Weyl group $W^\af = \mh Z R \rtimes W$ and in the
extended (affine) Weyl group $W^e = X \rtimes W$. As is well known, $W^\af$ is a 
Coxeter group, and the basis $\Pi$ of $R$ gives rise to a set $S^\af$ of simple
(affine) reflections. The length function $\ell$ of the Coxeter system 
$(W^\af ,S^\af )$ extends naturally to $W^e$. The elements of length zero form a
subgroup $\Omega \subset W^e$, and $W^e = W^\af \rtimes \Omega$. 

Let $q$ be a parameter function for $\mc R$, that is, a map $q : S^\af \to \mh C^\times$ 
such that $q(s) = q(s')$ if $s$ and $s'$ are conjugate in $W^e$. We also fix a square
root $q^{1/2} : S^\af \to \mh C^\times$.
The affine Hecke algebra $\mc H = \mc H (\mc R ,q)$ is the unique associative
complex algebra with basis $\{ N_w : w \in W^e \}$ and multiplication rules
\begin{equation}
\begin{array}{lll}
N_w \, N_{w'} = N_{w w'} & \mr{if} & \ell (w w') = \ell (w) + \ell (w') \,, \\
\big( N_s - q(s)^{1/2} \big) \big( N_s + q(s)^{-1/2} \big) = 0 & \mr{if} & s \in S^\af .
\end{array}
\end{equation}
In the literature one also finds this algebra defined in terms of the
elements $q(s)^{1/2} N_s$, in which case the multiplication can be described without
square roots. This explains why $q^{1/2}$ does not appear in the notation $\mc H (\mc R ,q)$.

In $X$ we have the positive cone
\[
X^+ := \{ x \in X : \inp{x}{\alpha^\vee} \geq 0 \; \forall \alpha \in \Pi \} \,.
\]
For $x \in X$ and $y,z \in X^+$ with $x = y - z$, we define $\theta_x := N_y N_z^{-1}$.
This is unambiguous, since $\ell$ is additive on $X^+$.
The span of the elements $\theta_x$ with $x \in X$ is a commutative subalgebra $\mc A$
of $\mc H$, isomorphic to $\mh C [X]$. We define a naive action of the group $W$ on
$\mc A$ by $w (\theta_x ) = \theta_{w (x)} $. Let $\mc H (W,q)$ be the finite 
dimensional Iwahori--Hecke algebra corresponding to the Weyl group $W$ and the 
parameter function $q |_S$.

\begin{thm}\label{thm:3.1}
\textup{(Bernstein presentation)}
\begin{description}
\item[a)] The multiplication in $\mc H$ induces isomorphisms of vector spaces
$\mc A \otimes \mc H (W,q) \to \mc H (\mc R ,q)$ and 
$\mc H (W,q) \otimes \mc A \to \mc H (\mc R ,q)$.
\item[b)] The center of $\mc H$ is $\mc A^W$, the invariants in $\mc A$ 
under the naive $W$-action.
\item[c)] For $f \in \mc A$ and $\alpha \in \Pi$ the following
Bernstein--Lusztig--Zelevinski relations hold:
\end{description}
\vspace{-3mm}
\[
f N_{s_\alpha} - N_{s_\alpha} s_\alpha (f) = \left\{ \!\! \begin{array}{ll}
\big( q (s_\alpha )^{1/2} - q(s_\alpha )^{-1/2} \big) (f - s_\alpha (f)) 
(\theta_0 - \theta_{-\alpha} )^{-1} & \alpha^\vee \notin 2Y \\
\big( q (s_\alpha )^{1/2} - q(s_\alpha )^{-1/2} + \! \big( q (\tilde s)^{1/2} - 
q( \tilde s )^{-1/2} \big) \theta_{-\alpha} \big) {\ds \frac{f - s_\alpha (f)}
{\theta_0 - \theta_{-2\alpha } }} & \alpha^\vee \in 2Y ,
\end{array} \right.
\]
\ \quad where $\tilde s \in S^\af$ is as in \textup{\cite[2.4]{Lus-Gr}}.
\end{thm}
\emph{Proof.}
See \cite[Section 3]{Lus-Gr}. $\qquad \Box$
\\[3mm]

The characters of $\mc A$ are elements of the complex torus 
$T = \mr{Hom}_{\mh Z}(X,\mh C^\times )$, which decomposes into a unitary and a
positive part:
\[
T = T_u \times T_{rs} = \mr{Hom}_{\mh Z}(X, S^1 ) \times \mr{Hom}_{\mh Z} (X,\mh R_{>0} ) \,.
\]
For any set $P \subset \Pi$ of simple roots we define
\[
\begin{array}{llllll}
X_P & = & X / (X \cap (P^\vee )^\perp ) & X^P & = & X / (X \cap \mh Q P ) \,, \\
Y_P & = & Y \cap Q P^\vee & Y^P & = & Y \cap P^\perp \,, \\
T_P & = & \mr{Hom}_{\mh Z} (X_P ,\mh C^\times ) & 
 T^P & = & \mr{Hom}_{\mh Z} (X^P ,\mh C^\times ) \,, \\
\mc R_P & = & (X_P ,R_P ,Y_P ,R_P^\vee ,P) & \mc R^P & = & (X,R_P ,Y,R_P^\vee ,P) \,, \\
\mc H_P & = & \mc H (\mc R_P ,q) & \mc H^P & = & \mc H (\mc R^P ,q) \,.
\end{array}
\]
The Lie algebras of $T_P$ and $T^P$ are $\mf t_P$ and $\mf t^P$, while the real forms 
$\mf a_P$ and $\mf a^P$ correspond to positive characters.

For $t \in T^P$ there is a surjective algebra homomorphism
\begin{equation}
\begin{aligned}
& \phi_t : \mc H^P \to \mc H_P \\
& \phi_t (N_w \theta_x ) = t(x) N_w \theta_{x_P} \,,
\end{aligned}
\end{equation}
where $x_P$ is the image of $x$ under the projection $X \to X_P$.

These constructions allow us to define parabolic induction.
For $t \in T^P$ and any $\mc H_P$-module $(\rho ,V_\rho)$ we put
\[
\pi (P,\rho ,t) = \mr{Ind}_{\mc H^P}^{\mc H} (\rho \circ \phi_t )
= \mc H \otimes_{\mc H^P} V_{\rho_t} \,.
\]
For any $\mc H$-module $(\pi ,V)$ and any $t \in T$ we have the $t$-weight space
\[
V_t := \{ v \in V : \pi (a) v = t(a) v \; \forall a \in \mc A \} \,.
\]
We say that $t$ is a weight of $V$ if $V_t \neq 0$. A finite dimensional $\mc H$-module
is called tempered if $|t(x)| \leq 1$ for all $x \in X^+$ and for all weights $t$.
If moreover $V$ is irreducible and $|t(x)| < 1$ for all $x \in X^+ \setminus \{0\}$ 
and for all weights $t$, then $V$ is said to belong to the discrete series.

There is a Langlands classification for irreducible modules of an affine Hecke algebra,
which is completely analogous to Theorem \ref{thm:2.1}. However, since it is more
awkward to write down, we refrain from doing so, and refer to \cite[Theorem 3.7]{Sol}.

Recall that the global dimension of $\mc H$ is the largest integer 
$d \in \mh Z_{\geq 0}$ such that the functor $\mr{Ext}_{\mc H}^d$ is not identically zero, 
or $\infty$ if no such number exists. We denote it by gl.\,dim$(\mc H)$.
It is known that gl.\,dim$(\mc H (\mc R ,q)) = \infty$ when the values of $q$ are certain 
roots of unity, but those are exceptional cases:

\begin{thm}\label{thm:3.5}
Suppose that 1 is the only root of unity in the subgroup of $\mh C^\times$ generated by
$\{ q(s)^{1/2} : s \in S^\af \}$. Then the global dimension of $\mc H$ equals rk$(X) = 
\dim_{\mh C} (\mf t )$.
\end{thm}
\emph{Proof.}
See \cite[Proposition 2.4]{OpSo1}. Although in \cite{OpSo1} $q$ is assumed to be positive,
the proof goes through as long as the finite dimensional algebra $\mc H (W,q)$
is semisimple. According to \cite[Theorem 3.9]{Gyo} this is the case under the 
indicated conditions on $q. \qquad \Box$
\\[3mm]

In the remainder of this section we assume that $q$ is positive, that is, 
$q(s)^{1/2} \in \mh R_{>0}$ for all $s \in S^\af$. 
Our affine Hecke algebra is equipped with an involution and a trace, namely,
for $x = \sum_{w \in W^e} x_w N_w \in \mc H$:
\[
x^* = {\ts \sum_{w \in W^e}} \overline{x_w} N_{w^{-1}} 
\qquad \mr{and} \qquad \tau (x) = x_e \,.
\]
Under the assumption that $q$ takes only positive values, * is an anti-automorphism,
and $\tau$ is positive. According to \cite[Proposition 1.12]{Opd-Sp} we have
\begin{equation}
\theta_x^* = N_{w_0} \theta_{-w_0 (x)} N_{w_0}^{-1} \,,
\end{equation}
where $w_0$ is the longest element of $W$. Every discrete series representation is 
unitary by \cite[Corollary 2.23]{Opd-Sp}. Moreover unitarity and temperedness are preserved
under unitary parabolic induction:

\begin{prop}\label{prop:3.2}
Assume that $q$ is positive and let $P \subset \Pi$.
\begin{description}
\item[a)] $\pi (P,\rho ,t)$ is unitary if $\rho$ is a unitary $\mc H_P$-representation and 
$t \in T^P_u$.
\item[b)] $\pi (P,\rho ,t)$ is tempered if and only if $\rho$ is a tempered
$\mc H_P$-representation and $t \in T^P_u$.
\end{description}
\end{prop}
\emph{Proof.}
The "if" statements are \cite[Propositions 4.19 and 4.20]{Opd-Sp}. The "only if" part of
b) can be proved just like Lemma \ref{lem:2.3} 
$\qquad \Box$ \\[3mm]

The bitrace $\inp{x}{y} := \tau (x^* y)$ gives $\mc H$ the structure of a Hilbert algebra.
This is the starting point for the harmonic analysis of $\mc H$ \cite{Opd-Sp,DeOp1}, 
which we prefer not to delve into here. We will need some of its deep results though.

Consider the finite group 
\[
K_P := T^P \cap T_P = T^P_u \cap T_{P,u} \,.
\]
For $k \in K_P$ and $w \in W(P,Q)$ there are algebra isomorphisms
\begin{equation}
\begin{aligned}
& \psi_w : \mc H^P \to \mc H^Q ,\\
& \psi_w (\theta_x N_w ) = \theta_{w (x)} N_{w u w^{-1}} \,, \\
& \psi_k : \mc H^P \to \mc H^P ,\\
& \psi_k (\theta_x N_u ) = k(x) \theta_x N_u \,.
\end{aligned}
\end{equation}
These maps descend to isomorphisms $\psi_w : \mc H_P \to \mc H_Q$ and
$\psi_k : \mc H_P \to \mc H_P$. The weights of the $\mc H_Q$-representation
\[
wk (\delta) := \delta \circ \psi_k^{-1} \circ \psi_w^{-1}
\]
are of the form $w (k^{-1} s) \in T_Q$, with $s \in T_Q$ a weight of $\delta$. 
Because $k \in T_u \,, wk (\delta )$ belongs to the discrete series of $\mc H_Q$. 

The space $\Xi$ of induction data for $\mc H$ consists of all triples $\xi = (P,\delta,t)$
where $P \subset \Pi \,, \delta$ is a discrete series representation of $\mc H_P$ and
$t \in T^P$. We call $\xi$ unitary, written $\xi \in \Xi_u$, if $t \in T_u^P$.
We write $\xi \cong \eta$ if $\eta = (P,\delta',t)$ with $\delta \cong \delta'$ as
$\mc H_P$-representations.

Let $\mc W$ be the finite groupoid, over the power set of $\Pi$, with
\[
\mc W_{PQ} = W(P,Q) \times K_P \,.
\] 
With the above we can define a groupoid action of $\mc W$ on $\Xi$ by
\begin{equation}
wk \cdot (P,\delta ,t) = (Q, wk (\delta ), w(kt)) \,.
\end{equation}
The algebra
\[
\prefix{_\mc F}{}{\mc H} := \mh C (T/W) \otimes_{Z (\mc H )} \mc H
\]
contains elements $\tau_w \; (w \in W)$ which satisfy the analogue of Proposition
\ref{prop:4.1} In fact Lusztig \cite[Section 5]{Lus-Gr} proved these results 
simultaneously for $\tau_w$ and $\tilde \tau_w$.
Let $\sigma$ be a discrete series representation of $\mc H^Q$ that is equivalent to
$wk (\delta)$, and let $I_\delta^{wk} : V_\delta \to V_\sigma$ be a unitary map such that
\begin{equation}\label{eq:3.2}
I_\delta^{wk} (\delta (\phi_t h) v ) = 
\sigma (\phi_{wt} \circ \psi_w \circ \psi_k (h)) (I_\delta^{wk} v) \,.
\end{equation}
Now we have a well-defined map
\begin{align*}
& I_{wk} : \prefix{_\mc F}{}{\mc H} \otimes_{\mc H^P} V_\delta \to 
 \prefix{_\mc F}{}{\mc H} \otimes_{\mc H^Q} V_\sigma , \\
& I_{wk} (h \otimes v) = h \tau^{-1}_w \otimes I_\delta^{wk} (v)
\end{align*}
We note that $I_{wk}$ can be constructed not just for $\delta$, but for any finite
dimensional $\mc H_P$-representation. Nevertheless for the next result we need that $\delta$
belongs to the discrete series and that $q$ is positive.

\begin{thm}\label{thm:3.3}
The map $I_{wk}$ defines an intertwining operator
\[
\pi (wk,P,\delta,t) : \pi (P,\delta,t) \to \pi (Q,\sigma ,w (kt) ) \,,
\]
which is rational as a function of $t \in T^P$. It is regular and invertible on an
analytically open neighborhood of $T_u^P$ in $T^P$. Moreover $\pi (wk,P,\delta,t)$ is
unitary if $t \in T_u^P$.
\end{thm}
\emph{Proof.}
See \cite[Theorem 4.33 and Corollary 4.34]{Opd-Sp}. For his intertwiners Opdam uses 
elements $\imath^o_w$ which are not quite the same as Lusztig's $\tau_w$. But these
approaches are equivalent, so the results from \cite{Opd-Sp} also hold in our setting.
$\qquad \Box$ \\[2mm]

In fact such operators yield all intertwiners between tempered parabolically induced modules:

\begin{thm}\label{thm:3.4}
For $\xi ,\eta \in \Xi_u$ the vector space $\mr{Hom}_{\mc H} (\pi (\xi) ,\pi (\eta))$
is spanned by $\{ \pi (w,\xi ) : w \in \mc W , w(\xi ) \cong \eta \}$.
\end{thm}
\emph{Proof.}
See \cite[Corollary 4.7]{DeOp1}. We remark that this result relies on a detailed
study of certain topological completions of $\mc H. \qquad \Box$
\vspace{4mm}

\section{Lusztig's reduction theorem} 
\label{sec:5}

In \cite{Lus-Gr} Lusztig established a strong connection between the representation
theories of affine Hecke algebras and graded Hecke algebras. We will use this to identify all 
intertwining operators between parabolically induced modules and to determine the
global dimension of a graded Hecke algebra. Let
\[
\mc R = (X,R,Y,R^\vee ,\Pi)
\]
be a based root datum and let
\[
\tilde{\mc R} = (X \otimes_{\mh Z} \mh R ,R, Y \otimes_{\mh Z} \mh R ,R^\vee ,\Pi )
\]
be the associated degenerate root datum. We endow $\mc R$ with a parameter function 
$q : S^\af \to {\mh C}^\times$ and $\tilde{\mc R}$ with the parameters
\begin{equation}\label{eq:5.1}
k_\alpha = \left\{ \begin{array}{lll}
\log \big( q(s_\alpha ) \big) & \mr{if} & \alpha^\vee \notin 2 Y \\
\log \big( q(s_\alpha ) q (\tilde s) \big) / 2 & \mr{if} & \alpha^\vee \in 2 Y,
\end{array}\right.
\end{equation}
where $\alpha \in \Pi$ and $\tilde s \in S^\af$ is as in \cite[2.4]{Lus-Gr}.
In case $q$ is not positive we fix a suitable branch of the logarithm in \eqref{eq:5.1}.
Every $k$ can be obtained in this way, in general even from several $X$ and $q$.

Let $\mc F^{me}(M)$ denote the algebra of meromorphic functions on a complex analytic 
variety $M$. The exponential map $\mf t \to T$ induces an algebra homomorphism
\begin{equation}
\begin{split}
& \Phi : \mc F^{me}(T)^W \otimes_{Z (\mc H )} \mc H \to 
 \mc F^{me}(\mf t )^W \otimes_{Z (\mh H )} \mh H \\
& \Phi \big( {\ts \sum_{w \in W}} f_w \tau_w \big) =
{\ts \sum_{w \in W}} (f_w \circ \exp ) \tilde \tau_w \hspace{2cm} f_w \in \mc F^{me}(T)
\end{split}
\end{equation}
For $\lambda \in \mf t$ let $Z_\lambda (\mh H) \subset Z (\mh H)$ be the maximal ideal of
functions vanishing at $W \lambda$. Let $\widehat{Z (\mh H )_\lambda}$ be the formal 
completion of $Z (\mh H)$ with respect to $Z_\lambda (\mh H)$, and define 
\begin{equation}
\hat{\mh H}_\lambda := \widehat{Z (\mh H )_\lambda} \otimes_{Z (\mh H)} \mh H \,.
\end{equation}
Recall that there is a natural bijection between finite dimensional
$\hat{\mh H}_\lambda$-modules and finite dimensional $\mh H$-modules whose generalized 
$S (\mf t^* )$-weights are all in $W \lambda$.

Similarly for $t \in T$ we have the maximal ideal $Z_t (\mc H) \subset Z (\mc H)$ and the
formal completions $\widehat{Z (\mc H )_t}$ and 
\begin{equation}
\hat{\mc H}_t := \widehat{Z (\mc H )_t} \otimes_{Z (\mc H)} \mc H \,.
\end{equation}
Finite dimensional $\hat{\mc H}_t$-modules correspond bijectively to finite dimensional
$\mc H$-modules with generalized $\mc A$-weights in $W t$. A slightly simplified version
of Lusztig's (second) reduction theorem \cite{Lus-Gr} states:

\begin{thm}\label{thm:5.1}
Suppose that 1 is the only root of unity in the subgroup of $\mh C^\times$ generated by
$\{ q(s)^{1/2} : s \in S^\af \}$, and let $\lambda \in \mf t$ be such that
\begin{equation}\label{eq:5.2}
\inp{\alpha}{\Im \lambda} \notin \pi \mh Z \setminus \{0\} \qquad \forall \alpha \in R \,.
\end{equation} 
Then the map $\Phi$ induces an algebra isomorphism
$\Phi_\lambda : \hat{\mc H}_{\exp (\lambda )} \to \hat{\mh H}_\lambda $.\\
This yields an equivalence $\Phi_\lambda^* = \Phi^*$ between the categories of:
\begin{itemize}
\item finite dimensional $\mh H$-modules whose $S (\mf t^* )$-weights are all in $W \lambda$,
\item finite dimensional $\mc H$-modules whose $\mc A$-weigths are all in $W \exp (\lambda )$.
\end{itemize}
\end{thm}
\emph{Proof.}
See \cite[Theorem 9.3]{Lus-Gr}. Our conditions on the $q$ replace the assumption \cite[9.1]{Lus-Gr}.
$\qquad \Box$
\\
\begin{cor}\label{cor:5.2}
Let $V$ be a finite dimensional $\hat{\mh H}_\lambda$-module, with $\lambda$ as in \eqref{eq:5.2}. 
The $\mc H$-module $\Phi_\lambda^* (V)$ is tempered if and only if $V$ is. 
Furthermore $\Phi_\lambda^* (V)$ belongs to the discrete series if and only if $V$ is a 
discrete series $\mh H$-module.
\end{cor}
\emph{Proof.}
These observations are made in \cite[(2.11)]{Slo}. We provide the (easy) proof anyway.
Let $\lambda_1 ,\ldots ,\lambda_d \in W \lambda$ be the $S(\mf t^*)$-weights of $V$. 
By construction the $\mc A$-weights of $\Phi_\lambda^* (V)$ are precisely 
$\exp (\lambda_1 ), \ldots, \exp (\lambda_d ) \in W (\exp \lambda)$.
Notice that 
\[
\exp (\Re \lambda_i ) = |\exp (\lambda_i )| \in T_{rs}
\]
and that the exponential map restricts to homeomorphisms $\mf a^- \to \{ t \in T_{rs} : 
t(x) \leq 1 \; \forall x \in X^+ \}$ and $\mf a^{--} \to \{ t \in T_{rs} : t(x) < 1 \; 
\forall x \in X^+ \setminus 0 \} \,. \qquad \Box$
\\[2mm]

Theorem \ref{thm:5.1} can also be used to determine the global dimension 
of a graded Hecke algebra. Although it is quite possible that this can be done without 
using affine Hecke algebras, the author has not succeeded in finding such a more 
elementary proof.

\begin{thm}\label{thm:5.4}
The global dimension of $\mh H$ equals $\dim_{\mh C} (\mf t^*)$.
\end{thm}
\emph{Proof.}
In view of the isomorphism \eqref{eq:1.3} we may assume that 
\begin{equation}\label{eq:5.4}
\mh Z \{ k_\alpha : \alpha \in R \} \cap \mh R i = \{ 0 \} .
\end{equation}
Let $q : S^\af \to \mh C^\times$ be a parameter function such that \eqref{eq:5.1} is satisfied.
Then \eqref{eq:5.4} assures that the subgroup of $\mh C^\times$ generated by $\{ q(s)^{1/2} :
s \in S^\af \}$ contains no roots of unity except 1, which will enable us to apply Theorem 
\ref{thm:3.5}.

For any $\lambda \in \mf t$ the localization functor $U \mapsto \hat{U}_\lambda := 
\hat{\mh H}_\lambda \otimes_{\mh H} U$ is exact and satisfies
\[
\mr{Hom}_{\hat{\mh H}_\lambda} (\hat{U}_\lambda , \hat{V}_\lambda ) \cong 
\widehat{Z (\mh H )_\lambda} \otimes_{Z (\mh H)} \mr{Hom}_{\mh H} (U,V) .
\]
for all $\mh H$-modules $U$ and $V$. Therefore
\[
\mr{Ext}^n_{\hat{\mh H}_\lambda} (\hat{U}_\lambda , \hat{V}_\lambda ) \cong 
\widehat{Z (\mh H )_\lambda} \otimes_{Z (\mh H)} \mr{Ext}^n_{\mh H} (U,V)
\]
for all $n \in \mh Z_{\geq 0}$. Every $\hat{\mh H}_\lambda$-module $M$ is of the form 
$\hat{U}_\lambda$ for some $\mh H$-module $U$ (namely $U = M$), so 
gl.\,dim$(\hat{\mh H}_\lambda) \leq \text{gl.\,dim}(\mh H)$. On the other hand the 
$Z(\mh H)$-module $\mr{Ext}^n_{\mh H} (U,V)$ is nonzero if and only if 
\[
\widehat{Z (\mh H )_\lambda} \otimes_{Z (\mh H)} \mr{Ext}^n_{\mh H} (U,V) \neq 0
\]
for at least one $\lambda \in \mf t$. We conclude that gl.\,dim$(\mh H) = 
\sup_{\lambda \in \mf t} \text{gl.\,dim} (\hat{\mh H}_\lambda)$. The same reasoning shows 
that gl.\,dim$(\mc H) = \sup_{t \in T} \text{gl.\,dim} (\hat{\mc H}_t)$.
Localizing \eqref{eq:1.3} yields an isomorphism
\[
\hat{\mh H}_\lambda = \widehat{\mh H (\tilde{\mc R},k)}_\lambda \cong 
\widehat{\mh H (\tilde{\mc R},zk)}_{z \lambda} .
\]
For some $z \in \mh R_{>0}$ the pair $(z \lambda ,zk)$ satisfies \eqref{eq:5.2} and 
\eqref{eq:5.4}. Then we can apply Theorems \ref{thm:3.5} and \ref{thm:5.1} to deduce that
\begin{equation*}
\text{gl.\,dim} (\hat{\mh H}_\lambda)  = \text{gl.\,dim} 
\big(\widehat{\mc H (\mc R, q^z)}_{\exp ( z\lambda)}\big) 
\leq \text{gl.\,dim}(\mc H (\mc R, q^z)) = \mr{rk}(X) = \dim_{\mh C}(\mf t^*).
\end{equation*}
Thus gl.\,dim$(\mh H) \leq \dim_{\mh C}(\mf t^*)$.\\
For the reverse inequality, let $\mh C_\lambda$ be the onedimensional $S(\mf t^*)$-module
with character $\lambda \in \mf t$ and consider the $\mh H$-module
$I_\lambda := \mr{Ind}_{S(\mf t^*)}^{\mh H} (\mh C_\lambda)$.
By Frobenius reciprocity
\begin{equation}\label{eq:5.5}
\mr{Ext}^n_{\mh H}(I_\lambda ,I_\lambda) \cong \mr{Ext}^n_{S(\mf t^*)} (\mh C_t ,I_\lambda)
= \bigoplus\nolimits_{w \in W} \mr{Ext}^n_{S(\mf t^*)} (\mh C_\lambda ,\mh C_{w \lambda}) .
\end{equation}
It is well-known that $\mr{Ext}_{S (\mf t^*)}^{\dim_{\mh C}(\mf t^*)}
(\mh C_\lambda ,\mh C_\mu)$ is nonzero if and only if $\lambda = \mu$. 
Hence gl.\,dim$(\mh H) \geq \dim_{\mh C}(\mf t^*).
\qquad \Box$ \\[3mm]

Now we set out to identify intertwining operators for parabolically induced $\mh H$-modules 
with those for $\mc H$-modules. We assume that $q$ is positive, and hence that $k$ is real.
Let $(\rho ,V_\rho )$ be a finite dimensional representation of $\mh H_P$, 
let $w \in W(P,Q)$ and let $\sigma$ be equivalent to $\rho \circ \psi_w^{-1}$. Since 
\[
\Phi_P \circ \psi_w^{-1} = \psi_w^{-1} \circ \Phi_Q \,,
\]
the $\mc H^Q$-representation $(\Phi_P^* \rho ) \circ \psi_w^{-1}$ is equivalent to 
$\Phi_Q^* (\sigma )$. We want to compare $I_\rho^w$ and $I^w_{\Phi_P^* (\rho )}$. 
Assume that $\mu$ is such that all weights of $\rho_\mu$ satisfy \eqref{eq:5.2}. 
For $h \in \mc H^P, t \in T^P$ and $h' = \Phi_P (\phi_{\exp (\mu )} h)$ we have by definition
\begin{align*}
I^w_{\Phi_P^* (\rho )} (\rho (h') v) & = I^w_{\Phi_P^* (\rho )} (\Phi_P^* (\rho ) ( 
\phi_{\exp (\mu )} (h)) v )\\
& = (\Phi_Q^* (\sigma )) (\phi_{w (\exp \mu )} \circ \psi_w (h)) (I^w_{\Phi_P^* (\rho )} v) \\
& = \sigma (\Phi_Q \circ \psi_w \circ \phi_{\exp (\mu )} (h)) (I^w_{\Phi_P^* (\rho )} v) \\
& = \sigma (\psi_w \circ \Phi_P \circ \phi_{\exp (\mu )} (h)) (I^w_{\Phi_P^* (\rho )} v) \\
& = \sigma (\psi_w h') (I^w_{\Phi_P^* (\rho )} v) \,.
\end{align*}
Hence $I^w_{\Phi_P^* (\rho )}$ satisfies the same intertwining property as $I^w_\rho$.  
Since we need $I^w_{\Phi_P^* (\rho )}$ to be unitary if $\rho$ is discrete series, 
we always define
\begin{equation}\label{eq:5.3}
I^w_\rho :=  I^w_{\Phi_P^* (\rho )} : V_\rho \to V_\sigma \,.
\end{equation}
\begin{prop}\label{prop:5.3}
Assume that $\lambda$ and all weights of $\rho_\mu$ satisfy \eqref{eq:5.2}.
\begin{description}
\item[a)] For any finite dimensional $\hat{\mh H}^P_\lambda$-module $V$, the map
\[
\Phi_\lambda \otimes \mr{Id}_V :  
\hat{\mc H}_{\exp (\lambda )} \otimes_{\hat{\mc H}^P_{\exp (\lambda )}} V \to 
\hat{\mh H}_\lambda \otimes_{\hat{\mh H}^P_\lambda} V 
\]
\!\! provides an isomorphism between the $\mc H$-modules 
$\mr{Ind}_{\mc H^P}^{\mc H} (\Phi_\lambda^{P*} V)$ 
and $\Phi_\lambda^* (\mr{Ind}_{\mh H^P}^{\mh H} V)$.
\item[b)] The following diagram commutes
\[
\begin{array}{lll}
\pi (P,\Phi_P^* (\rho ), \exp (\mu )) & \xrightarrow{\Phi \otimes \mr{Id}_{V_\rho}} &
 \pi (P,\rho, \mu ) \\
\downarrow \scs{\pi (w,P,\Phi_P^* (\rho ),\exp (\mu ))} & 
& \downarrow \scs{\pi (w,P,\rho ,\mu )} \\
\pi (Q,\Phi_Q^* (\sigma ),w(\exp \mu )) & \xrightarrow{\Phi \otimes \mr{Id}_{V_\sigma}} &
\pi (Q,\sigma ,w(\mu )) \,.
\end{array}
\]
\end{description}
\end{prop}
\emph{Proof.}
a) is a more concrete version of \cite[Theorem 6.2]{BaMo1}. By definition 
\[
\tilde \tau_w \in \mh C (\mf t ) \otimes_{\mh C [\mf t]} \mh H^P = 
\mh C (\mf t )^{W_P} \otimes_{Z (\mh H^P)} \mh H^P
\]
for all $w \in W_P$, and similarly 
\[
\tau_w \in \mh C (T ) \otimes_{\mh C [T]} \mc H^P = 
\mh C (T )^{W_P} \otimes_{Z (\mc H^P)} \mc H^P .
\]
Hence $\Phi_\lambda$ restricts to an isomorphism $\hat{\mc H}^P_{\exp (\lambda )} \to 
\hat{\mh H}^P_\lambda$. Since the multiplication in $\hat{\mh H}^P_\lambda$ induces a bijection
$\mh C [W^P] \otimes \hat{\mh H}^P_\lambda \to \hat{\mh H}_\lambda$, the multiplication in
$\hat{\mc H}_{\exp (\lambda )}$ provides a bijection 
\[
\Phi_\lambda^{-1} \big( \mh C [W^P] \big) \otimes \hat{\mc H}^P_{\exp (\lambda )} 
\to \hat{\mc H}_{\exp (\lambda )} \,.
\]
Consequently we can realize $\mr{Ind}_{\mc H^P}^{\mc H} (\Phi_\lambda^{P*} V)$ 
on the vector space $\Phi_\lambda^{-1} \big( \mh C [W^P] \big) \otimes V$. 
Now it is clear that the map from the proposition is a bijection, so we need to check that 
it is an $\mc H$-module homomorphism. For $h \in \hat{\mc H}_{\exp (\lambda )}$ we have
\begin{align*}
(\Phi_\lambda \otimes \mr{Id}_V ) (\mr{Ind}_{\mc H^P}^{\mc H} 
  (\Phi_\lambda^{P*} V)(h) (h' \otimes v)) 
& = (\Phi_\lambda \otimes \mr{Id}_V ) (h h'\otimes v) \\
& = \Phi_\lambda (h) \Phi_\lambda (h') \otimes v \\
& = \Phi_\lambda^* (\mr{Ind}_{\mh H^P}^{\mh H} V)(h) (\Phi_\lambda (h') \otimes v) \\
& = \Phi_\lambda^* (\mr{Ind}_{\mh H^P}^{\mh H} V)(h) 
  (\Phi_\lambda \otimes \mr{Id}_V ) (h' \otimes v) \,.
\end{align*}
b) On the larger space $\prefix{_\mc F}{}{\mc H} \otimes_{\mc H^P} V$ we have
\begin{align*}
\pi (w,P,\rho ,\mu ) (\Phi \otimes \mr{Id}_{V_\rho}) (h \otimes v) 
& = \pi (w,P,\rho ,\mu )(\Phi (h) \otimes v) \\
& = \Phi (h) \tilde \tau_w^{-1} \otimes I_\rho^w (v) \\
& = \Phi (h \tau_w^{-1}) \otimes I_\rho^w (v) \\
& = (\Phi \otimes \mr{Id}_{V_\sigma}) (h \tau_w^{-1} \otimes I^w_{\Phi_P^* (\rho )} (v)) \\
& = (\Phi \otimes \mr{Id}_{V_\sigma}) \pi (w,P, \Phi_P^* (\rho ),\exp (\mu )) (h \otimes v) \,.
\end{align*}
By assumption the image of $\mc H \otimes_{\mc H^P} V$ under these maps is 
$\mh H \otimes_{\mh H^P} V. \qquad \Box$
\vspace{4mm}

\section{Unitary representations}
\label{sec:6}

From now on we will assume that $k$ is real valued, i.e. that $k_\alpha \in \mh R$ for all
$\alpha \in R$. This assumption enables us to introduce a nice involution on $\mh H$, and
to speak of unitary representations. Let $w_0$ be the longest element of $W$.
Following Opdam \cite[p. 94]{Opd-Ha} we define
\begin{equation}
\begin{array}{lll@{\quad}l}
w^* & = & w^{-1} & w \in W , \\
x^* & = & w_0 \cdot \overline{-w_0 (x)} \cdot w_0 & x \in \mf t^* \,,
\end{array}
\end{equation}
where conjugation is meant with respect to the real form $\mf a^*$.

\begin{lem}\label{lem:6.1}
This * extends to a sesquilinear, anti-multiplicative involution on $\mh H$.
\end{lem}
\emph{Proof.}
It is clear that this is possible on $\mh C [W]$. Since $p \mapsto \overline{-w_0 (p)}$ 
and $p \mapsto w_0 p w_0$ are $\mh R$-linear automorphisms of the commutative algebra
$S (\mf t^* )$, * extends in the required fashion to $S(\mf t^* )$. Now we can define
\[
(w p)^* = p^* w^{-1} \quad \text{for} \quad p \in S(\mf t^* ) , w \in W ,
\]
and extend it to a sesquilinear bijection $\mh H \to \mh H$. To prove that this is an
anti-multiplicative involution we turn to the cross relation \eqref{eq:1.2}. 
It suffices to show that 
\begin{equation}\label{eq:6.1}
( k_\alpha \inp{x}{\alpha^\vee} + s_\alpha s_\alpha (x) - x s_\alpha )^* =
k_\alpha \overline{\inp{x}{\alpha^\vee}} + s_\alpha (x)^* s_\alpha - s_\alpha x^*
\end{equation}
is zero. We may assume that $x \in \mf a$, so that \eqref{eq:6.1} becomes
\[
k_\alpha \inp{x}{\alpha^\vee} - w_0 \cdot (w_0 s_\alpha) (x) \cdot w_0 \cdot s_\alpha + 
s_\alpha w_0 \cdot w_0 (x) \cdot w_0 \,.
\]
Conjugation with $w_0$ yields
\[
k_\alpha \inp{x}{\alpha^\vee} - (w_0 s_\alpha) (x) \cdot w_0 s_\alpha w_0 +
w_0 s_\alpha w_0 \cdot w_0 (x) \,.
\]
Let $y = w_0 (x) \in \mf a^*$ and let $\beta$ be the simple root $-w_0 (\alpha)$.
With \eqref{eq:1.2} we get
\begin{align*}
k_\alpha \inp{x}{\alpha^\vee} - s_\beta (y) s_\beta + s_\beta y & = 
k_\alpha \inp{x}{\alpha^\vee} + k_\beta \inp{y}{\beta^\vee} \\
& = k_\alpha \inp{x}{\alpha^\vee} + k_\beta \inp{w_0 (x)}{-w_0 (\alpha )^\vee} \\
& = (k_\alpha - k_\beta )\inp{x}{\alpha^\vee} \,.
\end{align*}
Since the automorphism $-w_0$ preserves the irreducible components of $R$, the roots
$\alpha$ and $\beta$ are conjugate in $W$. Hence $k_\alpha = k_\beta$, and 
\eqref{eq:6.1} is indeed zero. $\qquad \Box$
\\[2mm]

We note that
\[
x^* = \overline{-x} \quad \text{for} \quad x \in \mf t^{\Pi *} ,
\]
so that the $S(\mf t^{\Pi *} )$-weights of a unitary $\mh H$-module lie in $i \mf a^{\Pi *}$.

\begin{thm}\label{thm:6.2}
\begin{description}
\item[a)] The central character of a discrete series representation is real, i.e. lies in
$\mf a / W$.
\item[b)] There are only finitely many equivalence classes of discrete series representations.
\item[c)] Discrete representations are unitary.
\end{description}
\end{thm}
\emph{Proof.}
In Lemma 2.13 and Corollary 2.14 of \cite{Slo} Slooten proved a) and b), in a somewhat 
different way as we do below. We note that it is essential that all $k_\alpha$ are real. 

Let $(\delta ,V_\delta )$ be a
discrete series representation of $\mh H$ with central character $W \lambda \in \mf t / W$. 
By Corollary \ref{cor:5.2}  $\Phi_\lambda^* (\delta )$ is a discrete series representation 
of $\mc H (\mc R ,q)$, for a parameter function $q$ that satisfies \eqref{eq:5.1}. By 
\cite[Corollary 2.23 and Lemma 3.3]{Opd-Sp} $\Phi_\lambda^* (\delta )$ is unitary, and 
$\exp (\lambda) \in T$ is a residual point in the sense of \cite[Section 7.2]{Opd-Sp}.
This implies that $\lambda \in \mf t$ is a residual point. By \cite[Section 4]{HeOp} there
are only finitely many residual points in $\mf t$, and they all lie in $\mf a$. Since there
are only finitely many inequivalent irreducible $\mh H$-modules with a given central 
character, this proves a) and b).

Furthermore \cite[Theorem 3.10]{HeOp} tells us that $-\lambda \in W \lambda$, so in the 
terminology of \cite{BaMo1} $(\lambda ,k)$ is a real Hermitian point. Therefore we may invoke 
\cite[Theorem 5.7 and Corollary 5.8]{BaMo1}, which say that there exists an invertible 
Hermitian element $m = m^* \in \hat{\mh H}_\lambda$ such that 
\[
\Phi_\lambda (b^* ) = m \Phi_\lambda (b)^* m^{-1} \qquad b \in \hat{\mc H}_{\exp (\lambda )} \,.
\]
Moreover $m$ depends continuously on $(\lambda ,k)$, while for $k=0$ we have $m = 1$.
Hence $m$ is in fact a strictly positive element. Endow $V_\delta$ with a Hermitian inner 
product such that
\[
\delta (\Phi_\lambda (b) )^* = \delta (\Phi_\lambda (b^* )) 
\qquad \forall b \in \hat{\mc H}_{\exp (\lambda)} \,.
\]
For any $h \in \hat{\mh H}_\lambda$ we have
\[
\delta (h)^* = \delta (m h^* m^{-1}) = \delta (m) \delta (h^* ) \delta (m)^{-1} .
\]
In particular $\delta (m)^* = \delta (m)$, so $\delta (m) \in \mr{End}_{\mh C} (V_\delta )$
is again strictly positive. Let $\delta (m)^{1/2}$ be its unique positive square root in the
finite dimensional $C^*$-algebra $\mr{End}_{\mh C} (V_\delta )$. Now
\[
\delta (m)^{1/2} \delta (h^*) \delta (m)^{-1/2} = 
\delta (m)^{-1/2} \delta (h)^* \delta (m)^{1/2} =
\big( \delta (m)^{1/2} \delta (h) \delta (m)^{-1/2} \big)^* \,,
\]
so $\rho (h) := \delta (m)^{-1/2} \delta (h) \delta (m)^{1/2}$ is a unitary representation of
$\mh H$ on $V_\delta$. Since $\rho$ is clearly equivalent to $\delta$, we conclude that
$\delta$ is unitary as well. $\qquad \Box$.
\\[2mm]

Now we will investigate when unitarity is preserved under parabolic induction. Notice that 
this is not automatic, because the inclusion $\mh H^P \to \mh H$ does not always preserve 
the *. We define a Hermitian form $\inp{}{}_W$ on $\mh C [W]$ by declaring $W$ to be an 
orthonormal basis. 

\begin{prop}\label{prop:6.3}
Let $(\rho ,V_\rho )$ be a finite dimensional $\mh H_P$-module and $\lambda \in \mf t^P$.
The unitary dual of $\pi (P,\rho ,\lambda )$ is $\pi (P,\rho^* ,\overline{-\lambda} )$, where
$\rho^*$ is the unitary dual of $\rho$. The pairing between $\mh C [W^P] \otimes V_\rho$ and
$\mh C [W^P] \otimes V_\rho^*$ is given by
\[
\inp{w \otimes v}{w' \otimes v'} = \inp{w}{w'}_W \, \inp{v}{v'} \,.
\]
In particular $\pi (P,\rho ,\lambda )$ is unitary if $\rho$ is unitary and 
$\lambda \in i \mf a^P$. 
\end{prop}
\emph{Proof.}
According to \cite[Corollary 1.4]{BaMo2} the unitary dual of $\pi (P,\rho ,\lambda )$ is
$\mr{Ind}_{\mh H^P}^{\mh H} ((\rho_\lambda )^*)$, with respect to the indicated pairing.
Recall that $\mh H^P = S(\mf t^P )\otimes \mh H_P$, and that its involution satisfies $x^* = 
\overline{-x}$ for $x \in \mf t^P$. Hence the unitary dual of $\mh C_\lambda \otimes V_\rho$ 
is $\mh C_{\overline{-\lambda}} \otimes V_\rho^*$. In particular $\rho_\lambda $ is unitary
if $\rho$ is unitary and $\lambda \in i \mf a^P. \qquad \Box$
\\[2mm]

An induction datum for $\mh H$ is a triple $\xi = (P,\delta ,\lambda )$ such that
$P \subset \Pi \,, \lambda \in \mf t^P$ and $\delta$ belongs to the discrete series of
$\mh H_P$. We denote the space of such induction data by~$\tilde \Xi$. A second induction
datum $\eta$ is equivalent to $\xi$, written $\xi \cong \eta$, if $\eta = (P,\delta' ,\lambda)$
with $\delta' \cong \delta$ as $\mh H_P$-representations.
The subsets of unitary, respectively positive, induction data are defined as
\begin{equation}
\begin{array}{lll}
\tilde \Xi_u & = & 
  \big\{ (P,\delta ,\lambda) \in \tilde \Xi : \lambda \in i \mf a^P \big\} \,, \\
\tilde \Xi^+ & = & 
  \big\{ (P,\delta ,\lambda) \in \tilde \Xi : \Re (\lambda ) \in \mf a^{P+} \big\} \,.
\end{array}
\end{equation}
Notice that the (partially defined) action of $W$ on $\tilde \Xi$ preserves $\tilde \Xi_u$,
but not $\tilde \Xi^+$.
An obvious consequence of Theorem \ref{thm:6.2} and Proposition \ref{prop:6.3} is:

\begin{cor}\label{cor:6.4}
For any unitary induction datum $\xi \in \tilde \Xi_u$ the $\mh H$-module $\pi (\xi )$
is unitary and completely reducible.
\end{cor}

Slooten showed that unitary induction data are very useful for the classification of the 
tempered spectrum of $\mh H$:

\begin{thm}\label{thm:6.5}
For every irreducible tempered $\mh H$-module $V$ there exists a unitary induction
datum $\xi \in \tilde \Xi_u$ such that $V$ is equivalent to a direct summand of $\pi (\xi)$.
\end{thm}
\emph{Proof.}
See \cite[Section 2.2.5]{Slo}. $\qquad \Box$
\\[2mm]

\section{Classifying the intertwiners}
\label{sec:7}

Recall that $k$ is assumed to be real.
According to Theorem \ref{thm:3.4} the intertwiners corresponding to elements of $\mc W$
exhaust all homomorphisms between unitary parabolically induced $\mc H$-modules.  
It turns out that a similar, slightly simpler statement holds for graded Hecke algebras.

\begin{thm}\label{thm:7.1}
Let $\xi = (P,\delta ,\lambda ) ,\eta = (Q,\sigma ,\mu ) \in \tilde \Xi_u$. All the 
intertwining operators 
\[
\{ \pi (w,P,\delta ,\lambda ) : w \in W(P,Q) , w (\xi ) \cong \eta \}
\]
are regular and invertible at $\lambda$, and they span 
$\mr{Hom}_{\mh H} (\pi (\xi ), \pi (\eta))$.
\end{thm}
\emph{Proof.}
For the moment we assume that all weights of $\pi (\xi)$ and $\pi (\eta)$ satisfy 
\eqref{eq:5.2}. This holds for all $(\lambda ,\mu )$ in a dense open (with respect to 
the analytic topology) subset of $i \mf a^P \times i \mf a^Q$, and in particular on a 
suitable open neighborhood of $(0,0) \in i \mf a^P \times i \mf a^Q$. 

By Theorem \ref{thm:5.1} Lusztig's map $\Phi$ induces a bijection 
\[
\mr{Hom}_{\mh H} (\pi (\xi) ,\pi(\eta )) \to 
\mr{Hom}_{\mc H} (\Phi^* \pi (\xi ), \Phi^* \pi (\eta )) \,.\
\]
In view of Theorem \ref{thm:3.4} and Proposition \ref{prop:5.3}.a 
the right hand side is spanned by
\begin{multline}
\{ \pi (w k' ,P,\Phi_P^* (\delta ),\exp (\lambda )) : w \in W(P,Q), k' \in K_P , \\
w k' (\Phi_P^* (\delta )) \cong \Phi_Q^* (\sigma ) , w ( k' \exp \lambda ) = \exp (\mu ) \} \,.
\end{multline}
Moreover, according to Theorem \ref{thm:3.3} all these operators are regular and invertible.
By Theorem \ref{thm:6.2}.a the central characters of $\Phi_P^* (\delta )$ and 
$\Phi_Q^* (\sigma )$ are in $T_{rs} / W$. Since $K_P \subset T_u$, we can only get a contribution 
from $w k'$ if $k' = 1$. Now Proposition \ref{prop:5.3}.b completes the proof, under the above 
assumption on $\lambda$ and $\mu$.

For general $\lambda$ and $\mu$ we use a small trick. Consider the isomorphism\\
$m_z \colon \mh H (\tilde{\mc R}, zk) \to \mh H (\tilde{\mc R},k)$ from \eqref{eq:1.3}. 
One easily checks that 
\[
m_z^* \pi (P,\delta ,\lambda ) = \pi (P,m_z^* (\delta ),z \lambda )
\]
and that the following diagram commutes whenever the horizontal maps are well-defined:
\begin{equation}\label{eq:7.1}
\begin{array}{ccc} 
\mh H (\tilde{\mc R},zk) \otimes_{\mh H^P (\tilde{\mc R},zk)} V_\delta &
\xrightarrow{\pi (w,P,m_z^* (\delta ), z\lambda )} &
\mh H (\tilde{\mc R},zk) \otimes_{\mh H^Q (\tilde{\mc R},zk)} V_\sigma \\
\downarrow \scs{m_z \otimes \mr{Id}_{V_\delta}} & & 
  \downarrow \scs{m_z \otimes \mr{Id}_{V_\sigma}} \\
\mh H (\tilde{\mc R},k) \otimes_{\mh H^P (\tilde{\mc R},k)} V_\delta &
\xrightarrow{\pi (w,P,\delta ,\lambda )} &
\mh H (\tilde{\mc R},k) \otimes_{\mh H^Q (\tilde{\mc R},k)} V_\sigma 
\end{array}
\end{equation}
Now let $z > 0$ be positive. Then $m_z$ is not only an algebra isomorphism, it also preserves
the * and the real form $\mf a^*$ of $\mf t^*$. Take $z$ so small that all weights of
$m_z^* \pi (\xi )$ and $m_z^* \pi (\eta )$ satisfy \eqref{eq:5.2}. As we saw above, all the
intertwiners 
\[
\{ \pi (w,P,m_z^* (\delta ), z\lambda ) : w \in W(P,Q) , w (\xi ) \cong \eta \}
\]
are regular and invertible at $z \lambda \in i \mf a^P$, and they span 
$\mr{Hom}_{\mh H (\tilde{\mc R} ,zk)} (m_z^* \pi (\xi ), m_z^* \pi (\eta ))$. In view of 
\eqref{eq:7.1}, the same holds for the operators
\begin{equation}\label{eq:7.2}
\{ \pi (w,P,\delta ,\lambda ) : w \in W(P,Q) , w (\xi ) \cong \eta \} \,. \qquad \Box
\end{equation}
\vspace{2mm}

Several properties of these intertwiners are as yet unknown, but can be suspected 
from the analogy with affine Hecke algebras. 

In general the linear maps \eqref{eq:7.2} are linearly dependent.
To study this in detail, it is probably possible to develop the theory of R-groups for 
graded Hecke algebras, in analogy with the R-groups for reductive $p$-adic groups and 
affine Hecke algebras \cite{DeOp2}. 

As mentioned on page \pageref{eq:4.5}, $\pi (u, w(\xi )) \circ \pi (w,\xi )$ need not 
equal $\pi (uw ,\xi )$. By \eqref{eq:5.3} they can only differ by some scalar factor of 
absolute value one. Whether or not there always exists a clever choice of the $I^w_\delta$, 
which makes $w \mapsto \pi (w,\xi )$ multiplicative, is not known to the author.

In view of Theorem \ref{thm:3.3} it is not unlikely that $\pi (w,\xi )$ is unitary if
$\xi \in \tilde \Xi_u$. Yet this does not follow from Proposition \ref{prop:5.3}, since
$\Phi$ does not preserve the *.

For general induction data Theorem \ref{thm:7.1} fails, but fortunately it does extend to
positive induction data. We note that by \cite[Section 1.15]{Hum} every induction datum is
$W$-associate to a positive one. For $\xi = (P,\delta ,\lambda ) \in \tilde \Xi^+$ we write
\begin{equation}
\begin{array}{cll}
P(\xi ) & = & \{ \alpha \in \Pi : \inp{\alpha}{\Re (\lambda )} = 0 \} \,, \\
\xi_u & = & \big( P,\delta ,\lambda |_{\mf t^*_{P(\xi )}} \big) \,.
\end{array}
\end{equation}
Let $\pi^{P(\xi )}$ and $\pi_{P(\xi )}$ denote the induction functors for the graded 
Hecke algebras $\mh H^{P(\xi )}$ and $\mh H_{P(\xi )}$.

\begin{prop}\label{prop:7.2}
Let $\xi = (P,\delta ,\lambda ) \in \tilde \Xi^+$.
\begin{description}
\item[a)] The $\mh H^{P(\xi )}$-module $\pi^{P(\xi )}(\xi )$ is completely reducible,
and its restriction $\pi_{P(\xi )}(\xi_u )$ to $\mh H_{P(\xi )}$ is tempered and unitary.
\item[b)] Let $\mh C_\mu \otimes \rho$ be an irreducible constituent of 
$\pi^{P(\xi )}(\xi )$. Then $\mu = \lambda |_{\mf t^{P(\xi ) *}}$ and $(P(\xi ),\rho ,\mu )$ 
is a Langlands datum.
\item[c)] The irreducible quotients of $\pi (\xi )$ are precisely the modules 
$L(P(\xi ),\rho ,\mu )$ with $\rho$ and $\mu$ as in \textup{b)}. These modules are
tempered if and only if $\xi \in \tilde \Xi_u$.
\item[d)] Every irreducible $\mh H$-module can be obtained as in \textup{c)}.
\end{description}
\end{prop}
\emph{Proof.}
These results were inspired by the corresponding statements for affine Hecke algebras,
which were proved in unpublished work of Delorme and Opdam.\\
a) By definition $\lambda |_{\mf t^*_{P(\xi )}} \in i \mf a_{P(\xi )}$, so by Lemma
\ref{lem:2.3} and Corollary \ref{cor:6.4} $\pi_{P(\xi )}(\xi_u )$ is a tempered and
unitary $\mh H_{P(\xi )}$-module. Moreover $S(\mf t^{P(\xi ) *} )$ acts on 
$\pi^{P(\xi )}(\xi )$ by the character $\lambda |_{\mf t^{P(\xi ) *}}$, so 
$\pi^{P(\xi )}(\xi )$ is a completely reducible $\mh H^{P(\xi )}$-module.\\
b) Since $\xi \in \tilde \Xi^+$ we have
\[
\inp{\alpha}{\Re (\lambda )} > 0 \qquad \forall \alpha \in \Pi \setminus P(\xi ) \,,
\]
that is, $\mu = \lambda |_{\mf t^{P(\xi ) *}}$ has real part in $\mf a^{P(\xi )++}$.\\
c) By the transitivity of induction 
\[
\pi (\xi ) = \mr{Ind}_{\mh H^{P(\xi )}}^{\mh H} \pi^{P(\xi )} (\xi ) \,,
\]
so for the first statement we can apply Theorem \ref{thm:2.1}.a. If $\xi \in \tilde \Xi_u$, 
then all constituents of $\pi (\xi )$ tempered by Corollary \ref{cor:6.4}
On the other, if $\xi \notin \tilde \Xi_u$ then Theorem \ref{thm:2.1}.b tells us that
$L(P(\xi ),\rho ,\mu )$ cannot be tempered.\\
d) In view Theorem \ref{thm:2.1}.c it suffices to show that every irreducible tempered module
of a parabolic subalgebra of $\mh H$ appears as a direct summand of $\pi^{P(\xi )}(\xi )$,
for some $\xi \in \tilde \Xi^+$. But this is Theorem \ref{thm:6.5}. $\qquad \Box$
\\[2mm]

The representations $\pi (\xi )$ and $\pi (\eta )$ may have common irreducible 
constituents even if $\xi$ and $\eta$ are not $W$-equivalent in $\tilde \Xi^+$. 
This ambiguity disappears if we take only their irreducible quotients into account.

\begin{prop}\label{prop:7.3}
Let $\xi = (P,\delta ,\lambda ) , \eta = (Q,\sigma ,\mu ) \in \tilde \Xi^+$.
\begin{description}
\item[a)] The representations $\pi (\xi )$ and $\pi (\eta)$ have a common irreducible 
quotient if and only if there is a $w \in W(P,Q)$ with $w (\xi ) \cong \eta$.
\item[b)] If \textup{a)} applies, then $P(\xi ) = P(\eta )$ and the functor 
$\mr{Ind}_{\mh H^{P(\xi )}}^{\mh H}$ induces an isomorphism
\[
\hspace{-1cm}
\mr{Hom}_{\mh H_{P(\xi )}} (\pi_{P(\xi )}(\xi_u ), \pi_{P(\xi )}(\eta_u )) = 
\mr{Hom}_{\mh H^{P(\xi )}} (\pi^{P(\xi )}(\xi ), \pi^{P(\xi )}(\eta )) \cong
\mr{Hom}_{\mh H} (\pi (\xi ), \pi (\eta ))
\]
\item[c)] The operators 
\[
\{ \pi (w,\xi ) : w \in W(P,Q) , w (\xi ) \cong \eta \}
\]
are regular and invertible, and they span $\mr{Hom}_{\mh H} (\pi (\xi ),\pi (\eta))$. 
\end{description}
\end{prop}
\emph{Proof.}
a) Suppose that $\pi (\xi)$ and $\pi (\eta)$ have a common irreducible quotient.
By Proposition \ref{prop:7.2}.c and Theorem \ref{thm:2.1}.b we must have $P(\xi ) = P(\eta )$
and $\lambda |_{\mf t^{P(\xi ) *}} = \mu |_{\mf t^{P(\xi ) *}}$, while 
$\pi_{P(\xi )} (\xi_u )$ and $\pi_{P(\xi )}(\eta_u )$ must have a common irreducible
constituent. Applying Theorem \ref{thm:7.1} to $\mh H_{P(\xi )}$ we find a 
$w \in W_{P(\xi )}(P,Q)$ such that $w (\xi_u ) \cong \eta_u$. But $\xi_u$ (respectively 
$\eta_u$) differs only from $\xi$ (respectively $\eta$) by an element of $\mf t^{P(\xi )}$, 
so $w(\xi ) \cong \eta$ as well.

Conversely, suppose that $w \in W(P,Q)$ and $w (\xi ) \cong \eta$. Since $\Re (\lambda )$ 
and $\Re (\mu )$ are both in $\mf a^+$, they are equal, and fixed by $w$. From the definition 
we see that $P(\xi ) = P(\eta )$. Together with \cite[Proposition 1.15]{Hum} this shows that 
$w \in W_{P(\xi )}$ and $w (\xi_u ) = \eta_u$. Due to Theorem \ref{thm:7.1} 
$\pi_{P(\xi )}(\xi_u )$ and $\pi_{P(\xi )}(\eta_u )$ are isomorphic. To apply Proposition 
\ref{prop:7.2}.c we observe that, since $w \in W_{P(\xi )}$, 
\[
\lambda |_{\mf t^{P(\xi )}} = w(\lambda ) |_{\mf t^{P(\xi )}} = \mu |_{\mf t^{P(\xi )}} \,.
\]
b) From the above and Proposition \ref{prop:7.2}.a we see that the $\mh H^{P(\xi )}$-modules 
$\pi (\xi_u )$ and $\pi (\eta_u )$ are equivalent and completely reducible, and that 
$S(\mf t^{P(\xi ) *})$ acts on both by the character $\lambda |_{\mf t^{P(\xi ) *}}$. 
Hence we can apply Proposition \ref{prop:2.4}.a.\\ 
c) follows from b) and Theorem \ref{thm:7.1}. $\qquad \Box$
\\[2mm]

We remark that the maps $\lambda \mapsto \pi (w,\xi )$ can nevertheless have singularities, 
see page \pageref{eq:4.5}. These can even occur if $\xi \in \tilde \Xi_+$ but, according to 
Proposition \ref{prop:7.3}.c, not if both $\xi$ and $w(\xi )$ are positive.
\vspace{4mm}

\section{Extensions by diagram automorphisms}
\label{sec:8}

An automorphism $\gamma$ of the Dynkin diagram of the based root system $(R,\Pi )$ is a 
bijection $\Pi \to \Pi$ such that
\begin{equation}
\inp{\gamma (\alpha )}{\gamma (\beta )^\vee} = \inp{\alpha}{\beta^\vee} 
\qquad \forall \alpha ,\beta \in \Pi \,.
\end{equation}
Such a $\gamma$ naturally induces automorphisms of $R, R^\vee ,\mf t_\Pi ,\mf t^*_\Pi$ and 
$W$. Moreover we will assume that $\gamma$ acts on $\mf t$ and $\mf t^*$. If $\gamma$ and 
$\gamma'$ act in the same way on $R$ but differently on $\mf t^\Pi$, then we will, sloppily, 
regard them as different diagram automorphisms.

Let $\Gamma$ be a finite group of diagram automorphisms of $(R,\Pi)$. Groups like
\[
W' := \Gamma \ltimes W
\]
typically arise from larger Weyl groups as the isotropy groups of points in some torus, 
or as normalizers of some parabolic subgroup \cite{How}.
For the time being $k$ need not be real, but we do have to assume that 
$k_{\gamma (\alpha )} = k_\alpha \; \forall \alpha \in \Pi , \gamma \in \Gamma$.
Then $\Gamma$ acts on $\mh H$ by the algebra homomorphisms
\begin{equation}
\begin{split}
& \psi_\gamma : \mh H \to \mh H \,, \\
& \psi_\gamma (x s_\alpha ) = \gamma (x) s_{\gamma (\alpha )} 
  \qquad x \in \mf t^* , \alpha \in \Pi \,.
\end{split}
\end{equation}
In this section we will generalize Proposition \ref{prop:7.3} to the crossed product
\[
\mh H' := \Gamma \ltimes \mh H \,.
\]
We remark that algebras of this type play an important role in the classification of 
irreducible representations of affine Hecke algebras and reductive $p$-adic groups. 
See Lusztig's first reduction theorem \cite[Section 8]{Lus-Gr}. In the appendix we relate 
the representation theories of $\mh H$ and $\Gamma \ltimes \mh H$.

For any finite dimensional $\mh H_P$-module $(\rho ,V_\rho )$ and 
$\lambda \in \mf t^P$ we define the $\mh H'$-module
\begin{equation}\label{eq:8.4}
\pi' (P,\rho ,\lambda ) := \mr{Ind}_{\mh H}^{\mh H'} \pi (P,\rho ,\lambda ) = 
\mr{Ind}_{\mh H^P}^{\mh H'} \pi^P (P,\rho ,\lambda )
\end{equation}
For every $\gamma \in \Gamma$ and $P \subset \Pi$ we have algebra isomorphisms
\begin{equation}
\begin{split}
& \psi_\gamma : \mh H_P \to \mh H_{\gamma (P)} \,, \\
& \psi_\gamma : \mh H^P \to \mh H^{\gamma (P)} \,, \\
& \psi_\gamma (x s_\alpha ) = \gamma x s_\alpha \gamma^{-1} = \gamma (x) s_{\gamma (x)} 
\qquad x \in \mf t^* , \alpha \in P \,.
\end{split}
\end{equation}
In this situation $\gamma (\mf t_P^\pm ) = \mf t^\pm_{\gamma (P)}$, so we can define
\begin{equation}\label{eq:8.9}
\gamma (P, \rho ,\lambda ) = (\gamma (P), \rho \circ \psi_\gamma^{-1} ,\gamma (\lambda )) \,.
\end{equation}

\begin{lem}\label{lem:8.2}
Let $\gamma, P, \rho ,\lambda$ be as above and let $w \in W(P,Q)$.
\begin{description}
\item[a)] The $\mh H'$-modules $\pi' (P, \rho ,\lambda )$ and $\pi' (\gamma w (P, \rho ,\lambda ))$ 
have the same irreducible subquotients, with the same multiplicities.
\item[b)] The global dimension of $\mh H'$ equals $\dim_{\mh C}(\mf t^*)$.
\end{description}
\end{lem}
\emph{Proof.} a) Like in \eqref{eq:A.8} we see that 
$\pi' (\gamma w (P, \rho ,\lambda )) \cong \pi' (w (P, \rho ,\lambda ))$. By Lemma \ref{lem:4.4} 
the $\mh H$-modules $\pi (P, \rho ,\lambda )$ and $\pi (w (P, \rho ,\lambda ))$ have the same 
irreducible subquotients, counted with multiplicities. The functor $\mr{Ind}_{\mh H}^{\mh H'}$ 
preserves this property, because it is exact.\\
b) The calculation \eqref{eq:5.5}, with $\mh H$ replaced by $\mh H'$, shows that for 
$\lambda \in \mf t$
\[
\mr{Ext}^{\dim_{\mh C}(\mf t^*)}_{\mh H'} \big(\mr{Ind}_{S(\mf t^*)}^{\mh H'} (\mh C_\lambda),
\mr{Ind}_{S(\mf t^*)}^{\mh H'} (\mh C_\lambda) \big) \neq 0 ,
\]
so gl.\,dim$(\mh H') \geq \dim_{\mh C}(\mf t^*)$.

We claim that any $\mh H'$-module $U$ is isomorphic to a direct summand of 
$\mr{Ind}_{\mh H}^{\mh H'}(U)$. Consider the $\mh H'$-module homomorphism
\begin{align*}
& \phi : U \to \mr{Ind}_{\mh H}^{\mh H'}(U) , \\
& \phi (u) = |\Gamma|^{-1} \sum\nolimits_{\gamma \in \Gamma} 
\gamma \otimes \gamma^{-1} u .
\end{align*}
Cleary $\phi$ is injective, so $U \cong \phi (U)$ as $\mh H'$-modules. On the other hand there
is the multiplication map
\begin{align*}
& \mu : \mr{Ind}_{\mh H}^{\mh H'}(U) \to U ,\\
& \mu (h' \otimes u) = h' u .
\end{align*}
As $\mu \circ \phi = \mr{id}_U$, we find $\mr{Ind}_{\mh H}^{\mh H'}(U) = 
\phi (U) \oplus \ker (\mu)$. Hence, for every $\mh H'$-module~$V$
\[
\mr{Ext}_{\mh H'}^n (U,V) \cong \mr{Ext}_{\mh H'}^n (\phi(U),V) \subset
\mr{Ext}_{\mh H'}^n (\mr{Ind}_{\mh H}^{\mh H'}(U),V) \cong \mr{Ext}_{\mh H}^n (U,V). 
\]
For $n > \dim_{\mh C}(\mf t^*)$ the right hand side vanishes by Theorem \ref{thm:5.4},
so the left hand side is zero as well. We conclude that 
gl.\,dim$(\mh H') \leq \dim_{\mh C}(\mf t^*). \qquad \Box$
\\[2mm]

The elements of $\Gamma$ will give rise to new intertwiners between the 
$\pi' (P, \rho ,\lambda )$, in addition to the operators 
\[
\pi' (w,P, \rho ,\lambda ) :=  \mr{Ind}_{\mh H}^{\mh H'} \pi(w,P, \rho ,\lambda ) 
\qquad \mr{for} \; w \in W(P,Q) \,.
\]
Let $\sigma$ be a $\mh H_{\gamma (P)}$-representation which is equivalent to 
$\rho \circ \psi^{-1}_\gamma$, and let $I_\rho^\gamma : V_\rho \to V_\sigma$ be a linear
bijection such that
\begin{equation}\label{eq:8.7}
I_\rho^\gamma (\rho_\lambda (h) v) = \sigma_{\gamma (\lambda )}(\psi_\gamma h) 
(I_\rho^\gamma v) \qquad \forall h \in \mh H^P , v \in V_\rho \,.
\end{equation}
We define the intertwiner
\begin{equation}\label{eq:8.3}
\begin{split}
& \pi' (\gamma ,P,\rho ,\lambda ) : \pi' (P,\rho ,\lambda ) \to 
  \pi' (\gamma (P),\sigma ,\gamma (\lambda )) \,, \\
& \pi' (\gamma ,P,\rho ,\lambda ) (h \otimes v) = h \gamma^{-1} \otimes I_\rho^\gamma (v)
  \hspace{2cm} h \in \mh H', v \in V_\rho \,.
\end{split}
\end{equation}
Notice that this is simpler than \eqref{eq:4.5} and \eqref{eq:4.6}; $\pi' (\gamma ,\xi )$ 
is automatically well-defined and invertible, so there is no need to go to 
$\prefix{_\mc F}{}{\mh H}$ or $\prefix{_\mc F}{}{\mh H} \rtimes \Gamma$. We write 
\[
W' (P,Q) = \{ u \in W' : u (P) = Q \} \,.
\]
Since $\Gamma \cdot Q \subset R^+$, we have $W' (P,Q) \subset \Gamma W^P$.
For every $\gamma w \in W' (P,Q)$ we have a rational intertwiner
\begin{equation}
\pi' (\gamma w,P, \rho ,\lambda ) = 
 \pi' (\gamma ,w(P, \rho ,\lambda )) \circ \pi' (w,P, \rho ,\lambda )) 
\end{equation}
between $\pi (P, \rho ,\lambda )$ and $\pi (\gamma w (P, \rho ,\lambda ))$.

From now on we need to assume again that $k$ is real. Then \eqref{eq:8.9} defines an 
action of $\Gamma$ on $\tilde \Xi$, which preserves $\tilde \Xi^+$ and $\tilde \Xi_u$.

\begin{thm}\label{thm:8.1}
Let $\xi = (P,\delta ,\lambda) , \eta = (Q,\sigma ,\mu ) \in \tilde \Xi^+$.
\begin{description}
\item[a)] The representations $\pi' (\xi )$ and $\pi' (\eta )$ have a common irreducible 
quotient if and only if there is a $u \in W' (P,Q)$ with $u (\xi ) \cong \eta$.
\item[b)] The operators
\[
\{ \pi' (u,\xi ) : u \in W' (P,Q) , u (\xi ) \cong \eta \}
\]
are regular and invertible, and they span $\mr{Hom}_{\mh H'} (\pi' (\xi ), \pi' (\eta ))$.
\end{description}
\end{thm}
\emph{Proof.}
a) If $u = \gamma w$ and $u (\xi ) \cong \gamma' (\eta )$, then $w (\xi ) \cong \gamma^{-1} 
\gamma' (\xi ') \in \tilde \Xi^+$. So by Proposition \ref{prop:7.3}.c the intertwining 
operator $\pi'(u,\xi )$ is regular and invertible. In particular $\pi'(\xi )$ and 
$\pi' (\eta)$ are isomorphic $\mh H'$-representations, so they clearly have equivalent 
irreducible quotients.

Conversely, if $\pi'(\xi )$ and $\pi' (\eta )$ have a common irreducible quotient 
$\mh H'$-module, then they certainly have a common irreducible quotient $\mh H$-module. 
From \eqref{eq:A.10} we get an $\mh H$-module isomorphism
\begin{equation}\label{eq:8.1}
\pi'(\xi) \cong {\ts \bigoplus_{\gamma \in \Gamma}} \pi (\gamma \xi ) \,,
\end{equation}
and similary for $\eta$. From Proposition \ref{prop:7.3}.a and Lemma \ref{lem:A.2} we see that 
the $\mh H'$-modules $\pi' (\xi )$ and $\pi' (\eta)$ have a common irreducible quotient if and 
only if there exist $\gamma ,\gamma' \in \Gamma$ and $w \in W(\gamma' (P),\gamma (P))$ such that 
$w \gamma' (\xi) \cong \gamma (\eta)$. But this condition is equivalent to the existence of
$u = \gamma^{-1} w \gamma' \in W'(P,Q)$ with $u(\xi ) \cong \eta$.\\ 
b) The regularity and invertibility were already shown in the proof of part a).
By Frobenius reciprocity and \eqref{eq:8.1} we have
\begin{equation}\label{eq:.8.2}
\mr{Hom}_{\mh H'} (\pi' (\xi ), \pi' (\eta )) \cong \mr{Hom}_{\mh H} (\pi (\xi ), 
{\ts \bigoplus_{\gamma \in \Gamma}} \pi (\gamma \eta )) \,.
\end{equation}
By Proposition \ref{prop:7.3}.c the right hand side equals 
\[
{\ts \bigoplus_{\gamma \in \Gamma}} \mh C \{ \pi (w,\xi ) : w \in W(P,\gamma (Q)) , 
w(\xi ) \cong \gamma (\xi ') \} \,.
\]
Under the isomorphism \eqref{eq:8.1} $\mr{Ind}_{\mh H}^{\mh H'} \pi (w,\xi )$ 
corresponds to $\pi' (\gamma^{-1},w(\xi )) \circ \pi' (w,\xi ))$. Hence 
$\mr{Hom}_{\mh H'} (\pi' (\xi ),\pi' (\eta ))$ is spanned by the $\pi' (\gamma^{-1} w,\xi )$ 
with $w(\xi ) \cong \gamma (\eta )$, or equivalently $\gamma^{-1}w (\xi ) \cong \eta . 
\qquad \Box$ \\[2mm]

Now we can give a partial parametrization of irreducible $\mh H'$-modules,
in terms of our induction data.

\begin{thm}\label{thm:8.3}
Let $V$ be an irreducible $\mh H'$-module. There exists a unique association class
$W' (P,\delta ,\lambda ) \in \tilde \Xi / W'$ such that the following equivalent 
statements hold:
\begin{description}
\item[a)] $V$ is isomorphic to an irreducible quotient of $\pi' (\xi^+)$, for some
$\xi^+ \in \tilde \Xi^+ \cap W' (P,\delta ,\lambda )$,
\item[b)] $V$ is a constituent of $\pi' (P,\delta ,\lambda )$ and $\| cc_P (\delta ) \|$
is maximal for this property.
\end{description}
\end{thm}
\emph{Proof.}
By Proposition \ref{prop:7.2}.d and Theorem \ref{thm:A.1} there exists $\xi^+$ with
property a). By Theorem \ref{thm:8.1}.a the class $W' \xi^+ \in \tilde \Xi / W'$ is unique.\\
Pick $\xi = (P,\delta ,\lambda ) \in \tilde \Xi$ such that $V$ is a constituent of 
$\pi' (\xi )$ and $\| cc_P (\delta ) \|$ is maximal under this condition. 
By Lemma \ref{lem:8.2}.a we may assume that $\xi \in \tilde \Xi^+$. Suppose now that $V$
is not isomorphic to any quotient of $\pi' (\xi)$. As in Proposition \ref{prop:7.2}, let 
$\rho$ be an irreducible summand of the completely reducible $\mh H_{P (\xi)}$-module
$\pi_{P(\xi )} (\xi )$, such that $V$ is a constituent of 
\[
\mr{Ind}_{\mh H}^{\mh H'} \pi \big( P(\xi ), \rho ,\lambda |_{\mf t^{P(\xi) *}} \big) .
\]
By Theorem \ref{thm:6.2}.a and the definition of $P(\xi )$
\begin{equation}\label{eq:8.5}
\Re (cc_{P(\xi )} (\rho )) = W_{P(\xi )} cc_P (\delta ) .
\end{equation}
Choose an irreducible subquotient $\rho'$ of the $\mh H$-module 
$\pi \big( P(\xi ), \rho ,\lambda |_{\mf t^{P(\xi) *}} \big)$, such that $V$ is a
summand of $\mr{Ind}_{\mh H}^{\mh H'} (\rho' )$. With Theorem \ref{thm:2.1} we can associate
a Langlands datum $(Q,\sigma ,\mu )$ to $\rho'$, and by Proposition \ref{prop:2.4}.b
\begin{equation}\label{eq:8.6}
\| \Re (cc_{P (\xi )} (\rho )) \| < \| \Re (cc_Q (\sigma )) \| .
\end{equation}
Using Proposition \ref{prop:7.2}, find an induction datum $(Q' ,\sigma' ,\mu' ) \in \tilde \Xi^+$,
such that $\rho'$ is a quotient of $\pi (Q' ,\sigma' ,\mu' )$. Then 
\[
\Re (cc_Q (\sigma )) = W_Q cc_{Q'} (\sigma' ) ,
\]
which together with \eqref{eq:8.5} and \eqref{eq:8.6} yields 
\[
\| cc_P (\delta ) \| < \| cc_{Q'} (\sigma' ) \| .
\]
But by construction $V$ is a constituent of $\pi' (Q',\sigma' ,\mu')$, so this is a 
contradiction. We conclude that $V$ has to be a quotient of $\pi' (P,\delta ,\lambda )$. 

Hence the association class $W' (P,\delta ,\lambda )$ satisfies not only b) but also a), 
which at the same time shows that it is unique. In particular the conditions a) and b)
turn out to be equivalent. $\qquad \Box$
\vspace{4mm}

\section{The space of module homomorphisms}
\label{sec:9}

How does the algebra $\mr{End}_{\mh H'} (\pi' (P,\delta ,\lambda ))$ change if we vary 
$\lambda \in \mf t^P$? Proposition \ref{prop:2.4} strongly suggests that it is rigid if
we insist that $\Re (\lambda ) \in \mf a^{P++}$. Although the intertwiners 
$\pi' (w,P,\delta ,\lambda )$ are by no means constant, their span tends to be stable.
Of course $\mr{End}_{\mh H'} (\pi' (P,\delta ,\lambda ))$ might jump if the isotropy
group of $\lambda$ in $\Gamma \ltimes W$ becomes larger, but then it should only grow.
This and even more turns out to hold true:

\begin{thm}\label{thm:9.1}
Let $\xi = (P,\delta ,\lambda ) ,\eta = (Q,\sigma ,\mu ) \in \tilde \Xi^+$. Suppose that 
$\lambda' \in \mf t^P , \mu' \in \mf t^Q$ and $u (\lambda' ) = \mu'$ for all $u \in W'(P,Q)$ 
with $u (\lambda ) = \mu$. Then, as subsets of $\mr{End}_{\mh C}(\mh C [\Gamma W^P] 
\otimes V_\delta , \mh C [\Gamma W^Q] \otimes V_\sigma )$, we have
\[
\mr{Hom}_{\mh H'} (\pi' (\xi ),\pi' (\eta )) \subset \mr{Hom}_{\mh H'} 
(\pi' (P,\delta ,\lambda' ), \pi' (Q,\sigma ,\mu' )) \,.
\]
\end{thm}
\emph{Proof.}
In view of Theorem \ref{thm:8.1}.c it suffices to consider
\[
f = \pi (u,\xi ) \in \mr{Hom}_{\mh H'} (\pi' (\xi ),\pi' (\eta )) \,.
\]
Our conditions imply that $u (P,\delta ,\lambda' ) \cong (Q,\sigma ,\mu' )$. 

As $\Gamma \ltimes W$-representations we have 
\[
\begin{array}{lll}
\pi' (\xi ) & = & \mr{Ind}_{W_P}^{\Gamma \ltimes W} \delta \,, \\
\pi' (\eta ) & = & \mr{Ind}_{W_Q}^{\Gamma \ltimes W} \sigma \,,
\end{array}
\]
and the operators $\pi' (\xi ) (\gamma w)$ and $\pi' (\eta ) (\gamma w)$ 
do not depend on $\lambda$ and $\mu$. Hence 
\begin{equation}\label{eq:9.4}
f \in \mr{Hom}_{\Gamma \ltimes W} (\pi' (P,\delta ,\lambda' ), \pi' (Q,\sigma ,\mu' )) 
\quad \forall \lambda' \in \mf t^P ,\mu' \in \mf t^Q \,.
\end{equation}
Since $\Gamma \ltimes W$ and $\mf t^*$ generate $\mh H'$, it suffices to check that 
\begin{equation}\label{eq:9.1}
f (\pi' (P,\delta ,\lambda' )(x) \: w \otimes v) = \pi' (Q,\sigma ,\mu' )(x) \, f(w \otimes v) 
\quad \forall x \in \mf t^* , v \in V_\delta , w \in \Gamma W^P \,,
\end{equation}
for all $\lambda' ,\mu'$ as in the theorem. Moreover we claim that we may restrict to $w=1$. 
Indeed, if we know \eqref{eq:9.1} for $1 \otimes v$, then for any $w \in \Gamma W^P$ we get
\begin{multline}
f (\pi' (P,\delta ,\lambda' )(x) \, w \otimes v) = 
 f (\pi' (P,\delta ,\lambda' )(xw) \, 1 \otimes v) = \\
\pi' (Q,\sigma  ,\mu' )(xw) \, f(1 \otimes v) = \pi' (Q,\sigma  ,\mu' )(x) \, w f(1\otimes v) 
\end{multline}
and from the explicit formulas \eqref{eq:4.5} and \eqref{eq:8.3} we see that
\[
w f (1 \otimes v) = f (w \otimes v) \,.
\]
Let $w_1 , \ldots ,w_{|W^Q |}$ be the elements of $W^Q$, listed in a length-increasing way. 
Write $\Gamma = \{ \gamma_1 ,\ldots ,\gamma_{|\Gamma |} \}$ and let $v_1 , \ldots ,v_{\dim V}$ 
be a basis of $V_\sigma$. Then 
\[
\gamma_1 w_1 \otimes v_1 , \ldots , \gamma_{|\Gamma |} w_1 \otimes v_{\dim V} , \gamma_1 w_2 
\otimes v_1 ,\ldots , \gamma_{|\Gamma |} w_{|W^Q |} \otimes v_{\dim V} 
\] 
is a basis of $\mh H' \otimes_{\mh H^Q} (\mh C_{\mu} \otimes V_\sigma )$. The important thing 
here is that the elements of $W^Q$ appear from short to long. 

Let $s_1 \cdots s_r = w_m \in W^Q$ be a reduced expression, with $s_i = s_{\alpha_i} \in S$. 
By a repeated application of the cross relation \eqref{eq:1.2} we find for any 
$x \in \mf t^* ,\gamma \in \Gamma$:
\begin{multline}
x \gamma w_m =  \gamma \gamma^{-1}(x) w_m = 
(\gamma w_m )^{-1}(x) \, + \\ 
\gamma {\ts \sum_{i=1}^r} k_{\alpha_i} \inp{\gamma^{-1}(x)}{
s_1 \cdots s_{i-1} s_{i+1} \cdots s_r (\alpha_i^\vee )} s_1 \cdots s_{i-1} s_{i+1} \cdots s_r \,.
\end{multline}
If $(\gamma w_m )^{-1}(x) = (\gamma w_m )^{-1}(x)_Q + (\gamma w_m )^{-1}(x)^Q \in 
\mf t_Q \oplus \mf t^Q$, then
\begin{equation*}
\begin{array}{lll}
\sigma_\mu ((\gamma w_m )^{-1}(x)) v' & = & 
\! \inp{(\gamma w_m )^{-1}(x)^Q}{\mu} v' + \sigma ((\gamma w_m )^{-1}(x)_Q) v' \,, \\
\pi' (\eta )(x) (\gamma w_m (x) \otimes v') \!\!\! & = & 
 \!\gamma w_m \otimes \inp{(\gamma w_m )^{-1}(x)^Q \!}{\!\mu} 
v' + \gamma w_m \otimes \sigma ((\gamma w_m )^{-1}(x)_Q) v'  \\
\multicolumn{3}{c}{+ \, \sum_{i=1}^r k_{\alpha_i} \inp{x}{\gamma s_1 \cdots s_{i-1} 
s_{i+1} \cdots s_r (\alpha_i^\vee )} \gamma s_1 \cdots s_{i-1} s_{i+1} \cdots s_r \otimes v' \,.}
\end{array}
\end{equation*}
By definition all elements $s_1 \cdots s_{i-1} s_{i+1} \cdots s_r$ can be written as $w_j u$
with $w \in W_P$ and $j < m$. Hence we can express the above as
\begin{align}
\nonumber & \pi' (\eta )(x) (\gamma w_m (x) \otimes v' ) \:=\: 
\gamma w_m \otimes \sigma ((\gamma w_m )^{-1}(x)_Q) v' \: + \: 
 M(x,\mu ) (\gamma w_m (x) \otimes v') \,, \\
\label{eq:9.2} & M(x,\mu ) \:=\: \begin{pmatrix}
\inp{x}{\gamma_1 w_1 (\mu )} & & * & & * \\
 & \ddots & & & \\
 0 & & \inp{x}{\gamma_1 w_2 (\mu )} & & * \\
 & & & \ddots & \\
 0 & & 0 & & \inp{x}{\gamma_{|\Gamma |} w_{|W^Q |} (\mu )}
\end{pmatrix} \,,
\end{align}
where the stars denote expressions of the form $\inp{x}{\kappa}$ with $\kappa$ independent 
of $\mu$. From this matrix we see that all generalized weights of $\pi' (\eta )$ are of 
the form $\gamma w_m (\mu + \nu )$, with $\nu $ a $S(\mf t_W^* )$-weight of $\sigma$.

Now we will first finish the proof under the assumption $\lambda ,\mu \in i \mf a$, 
and then deal with general $\lambda$ and $\mu$.

We may assume that $v \in V_\delta$ is a generalized $S(\mf t^* )$-weight vector, 
with weight $\lambda + \nu'$. Because $1 \otimes v$ is cyclic for $\pi' (\xi ) \,, 
\pi' (u,\xi )(1 \otimes v) \neq 0$. Hence $\lambda + \nu' = \gamma w_m (\mu + \nu )$ 
for some $\gamma ,w_m ,\nu$ as above. In particular
\begin{equation}\label{eq:9.6}
\nu' = \Re (\lambda + \nu' ) = \Re (\gamma w_m (\mu + \nu )) = \gamma w_m (\nu ) \,.
\end{equation}
We note that $\nu \in \mf a_Q^{--} , \nu' \in \mf a_P^{--}$ and 
$\gamma w_m (Q) \subset R^+$. Therefore $\gamma w_m (Q) \subset R_P^+$. 
Because $|Q| = |P|$, it follows that actually $\gamma w_m \in W'(Q,P)$. This shows that
\begin{equation}
f (1 \otimes V_\delta ) \subset \mh C [W'(Q,P)] \otimes V_\sigma \,.
\end{equation} 
Furthermore $\gamma w_m (\mf t_Q ) = \mf t_P$, so for $x \in \mf t_P^*$ we have 
$(\gamma w_m )^{-1}(x)^Q = 0$. Hence $\inp{x}{\gamma w_m (\mu' + \nu )}$ and \
$\pi' (Q,\sigma ,\mu')(x) \, f(1 \otimes v)$ do not depend on $\mu'$. We conclude that 
\begin{equation}\label{eq:9.5}
\mr{Hom}_{S(\mf t_P^* )} (\mh C_\lambda \otimes V_\delta , \pi' (Q,\sigma ,\mu ')) =
\mr{Hom}_{S(\mf t_P^* )} (V_\delta , \mh C [W'(Q,P)] \otimes V_\sigma )
\end{equation}
is independent of $\mu'$.

For $x \in \mf t^{P*}$ we have $(\gamma w_m )^{-1}(x) \in \mf t^{Q*}$, 
so by \eqref{eq:9.2} and the definition of $\delta_\lambda$:
\begin{equation}\label{eq:9.3}
\begin{array}{lll}
\pi' (\eta )(x) \, f(v) & = & M(x,\mu ) f(v) \,, \\
\pi' (\xi )(x) v & = & \delta_\lambda (x) v = \inp{x}{\lambda} v \,.
\end{array}
\end{equation}
Therefore $f(1 \otimes v)$ lies in the $\inp{x}{\lambda}$-eigenspace of $M(x,\mu )$. 
The matrix of 
\[
x \mapsto M(x,\mu ,\lambda ) := M(x,\mu ) - \inp{x}{\lambda}
\] 
has diagonal entries $\gamma w (\mu ) - \lambda$, while its off-diagonal entries do 
not depend on $\lambda$ and $\mu$. We are interested in
\[
K(\mu ,\lambda ) := {\ts \bigcap_{x \in \mf t^{P*}}} \ker M(x,\mu ,\lambda ) \,.
\]
By Theorem \ref{thm:8.1}.c $\pi'(u,P,\delta ,\lambda )$ is a rational function of $\lambda$, 
and it is regular $\forall \lambda \in \mf t^{P+}$. Consequently $\pi'(u,P,\delta ,\lambda )(x) 
(1 \otimes v) \in K(u(\lambda ),\lambda )$ for all $\lambda \in \mf t^{P+}$. 
There is a minimal $K_0 \subset \mh C [\Gamma W^Q ] \otimes V_\sigma$, such that 
$K(u(\lambda ),\lambda ) = K_0$ for all $\lambda$ in a Zariski-open subset of $\mf t^P$. 
By continuity $\pi'(u,P,\delta ,\lambda ) (1 \otimes v) \in K_0$ for all $\lambda \in \mf t^{P*}$, 
and in particular $f (1 \otimes v) \in K_0$. Now our initial conditions on $\lambda'$ and $\mu'$ 
assure that
\begin{equation*}
f (\pi' (P,\delta ,\lambda' )(x) (1 \otimes v) = f (\inp{x}{\lambda'} \otimes v) = 
\inp{x}{\gamma w_m (\mu' )} f(1 \otimes v) = \pi' (Q,\sigma ,\mu' ) f(1 \otimes v) \,.
\end{equation*}
for all $x \in \mf t^{P*}$. Together with \eqref{eq:9.4} and \eqref{eq:9.5} this implies that 
\[
f : \mh C [\Gamma W^P] \otimes V_\delta \to \mh C [\Gamma W^Q] \otimes V_\sigma
\]
is in $\mr{Hom}_{\mh H'} (\pi' (P,\delta ,\lambda' ) ,\pi' (Q,\sigma ,\mu' ))$ for all eligible
$(\lambda' ,\mu')$. We note that it may not equal $\pi' (u,P,\delta ,\lambda' )$ everywhere.

As promised, we will now discuss general $\lambda \in \mf t^{P*}$ and $\mu \in \mf t^{Q*}$.
If $u = \gamma w \in \Gamma W^P$, then by definition 
$\pi' (u,\xi) = \pi'(\gamma, w (\xi)) \circ \pi' (w,\xi )$. Clearly
\[
\pi' (\gamma,w (\xi )) \in \mr{End}_{\mh C} \big( \mh C [\Gamma W^{w(P)}] \otimes V_\delta \big)
\]
is independent of $\lambda$. For any $r > 0$ we have 
\[
(P,\delta ,r \lambda )_u = \big( P,\delta ,r \lambda |_{\mf t^*_{P(\xi )}} \big) 
\in \tilde \Xi_{P(\xi ),u} \,,
\]
so by the above 
\[
\pi_{P(\xi )} (w,\xi ) \in \mr{Hom}_{\mh H_{P(\xi )}} (\pi_{P(\xi )} (P,\delta ,r \lambda )_u ,
\pi_{P(\xi )} w(P,\delta ,r \lambda )_u ) \,.
\]
By Proposition \ref{prop:7.3}.b we have
\[
\pi' (w,\xi ) = \mr{Ind}_{\mh H_{P(\xi )}}^{\mh H'} \pi_{P(\xi )} (w,\xi_u ) \,,
\]
so $\pi'(w,\xi ) \in \mr{Hom}_{\mh H'} (\pi' (P,\delta ,r \lambda ) 
,\pi' w(P,\delta ,r \lambda ))$ for all $r > 0$. Now we can follow the above proof up to 
\eqref{eq:9.6}, which we have to replace by
\[
r \Re (\lambda ) + \nu'= \Re (r \lambda + \nu' ) = 
\Re (\gamma w_m (r \mu + \nu )) = \gamma w_m (r \Re (\mu ) + \nu) \,.
\]
Since this holds for all $r > 0$, we conclude as before that $\nu' = \gamma w_m (\nu )$.
The rest of the proof goes through without additional problems. $\qquad \Box$
\\[2mm]

A particular case of Theorem \ref{thm:9.1} appears for every $\mh H'$-module decomposition
$\pi' (\xi ) = V_1 \oplus V_2$. The projections $p_i : \mh C [\Gamma W^P ] \otimes V_\delta 
\to V_i$ are in $\mr{End}_{\mh H'} (\pi' (\xi ))$, so they are also in 
$\mr{End}_{\mh H'} (\pi' (P,\delta ,\lambda' ))$
for all $\lambda' = \mu'$ as in Theorem \ref{thm:9.1}. This shows that 
$\pi' (P,\delta ,\lambda' ) = V_1 \oplus V_2$ is also a decomposition of $\mh H'$-modules. 
In particular this holds for $\lambda' = r \lambda$ with $r \geq 0$, from which we conclude 
that the maximal reducibility of $\pi' (P,\delta, \lambda' )$ occurs at $\lambda' = 0$.
\vspace{4mm}

\section{Geometric description of the spectrum}
\label{sec:10}

By the spectrum of an algebra $A$ we mean the collection Irr$(A)$ of equivalence 
classes of irreducible $A$-representations, endowed with the Jacobson topology. 
This topology generalizes the Zariski-topology for commutative rings, its closed
subsets are by definition of the form
\[
V(S) := \{ \pi \in \mr{Irr}(A) : \pi (s) = 0 \; \forall s \in S \} ,
\]
for any subset $S \subset A$.

In this section we look at the spectrum Irr($\mh H'$) of $\mh H' = \Gamma \ltimes \mh H$ 
from a geometric point of view. For every discrete series representation $\delta$ of a 
parabolic subalgebra $\mh H_P$ we get a series of $\mh H'$-modules 
$\pi' (P,\delta ,\lambda )$, parametrized by $\lambda \in \mf t^P$. The finite group
\[
W'_\delta := \{ u \in \Gamma \ltimes W : u (P) = P, u (\delta ) \cong \delta \}
\]
acts linearly on $\mf t^P$, and $\lambda$'s in the same $W'_\delta$-orbit yield modules 
with the same irreducible subquotients. Theorem \ref{thm:8.3} allows us to associate to 
every irreducible $\mh H'$-module $\rho$ an induction datum $\xi^+ (\rho) \in \tilde \Xi^+$, 
unique up to $W'$-equivalence, such that $\rho$ is a quotient of $\pi' (\xi^+ (\rho ))$. 
For any subset $U \subset \mf t^P$ we put
\[
\mr{Irr}_{(P,\delta ,U)} (\mh H' ) = \{ \rho \in \mr{Irr}(\mh H') : 
W' \xi^+ (\rho ) \cap (P,\delta ,U) \neq \emptyset \} .
\]
For $U = \mf t^P$ or $U = \{ \lambda \}$ we abbreviate this to 
$\mr{Irr}_{P,\delta} (\mh H' )$ or $\mr{Irr}_{(P,\delta ,\lambda)} (\mh H' )$.
The closures of the sets $\mr{Irr}_{P,\delta} (\mh H' )$ define a filtration of Irr$(\mh H)$,
which is the analogue of the stratification of the smooth dual of a reductive $p$-adic group 
from \cite[Lemma 2.17]{Sol3}.

\begin{prop}\label{prop:10.1}
Suppose that $U \subset \mf t^P$ satisfies
\begin{itemize}
\item $U \to \{-1,0,1\} : \lambda \mapsto \mr{sign} \inp{\Re (\lambda )}{\alpha}$
is constant for all $\alpha \in \Pi \setminus P$,
\item every $\lambda \in U$ has the same stabilizer in $W'_\delta$.
\end{itemize}
Then $\mr{Irr}_{(P,\delta ,U)} (\mh H' )$ is homeomorphic to $U / W'_\delta \times 
\mr{Irr}_{(P,\delta ,\lambda_0 )} (\mh H' )$ for any $\lambda_0 \in U$.
\end{prop}
\emph{Proof.}
The first assumption implies that there exists $w \in W$ such that $w (P,\delta ,U) 
\subset \tilde \Xi^+$. Hence we may assume without loss of generality that 
$U \subset \mf a^{P+} \oplus i \mf a^P$. Write 
\[
Q = \{ \alpha \in \Pi : \inp{\Re (\lambda )}{\alpha } = 0 \; \forall \lambda \in U \} .
\]
If $\lambda ,\lambda' \in U$ and $w \in W_Q$, then
\[
w \big( \lambda |_{\mf t_Q^*} \big) =  \lambda |_{\mf t_Q^*} \; \Longleftrightarrow \;
w ( \lambda ) = \lambda \; \Longleftrightarrow \; 
w ( \lambda' ) = \lambda' \; \Longleftrightarrow \; 
w \big( \lambda' |_{\mf t_Q^*} \big) = \lambda' |_{\mf t_Q^*} .
\]
By Corollary \ref{cor:6.4} the $\mh H_Q$-module 
$\pi_Q \big( P,\delta ,\lambda |_{\mf t_Q^*} \big)$ is completely reducible, 
and by Theorem \ref{thm:9.1} 
\[
E := \mr{End}_{\mh H_Q} \big( \pi_Q \big( P,\delta ,\lambda |_{\mf t_Q^*} \big) \big)
\]
does not depend on $\lambda \in U$. Consequently $E$ is a semisimple algebra and every
$\mh H_Q$-submodule of $\pi_Q \big( P,\delta ,\lambda |_{\mf t_Q^*} \big)$ is of the
form im$(e)$ for some idempotent $e \in E$. Such a module is irreducible if and only if
$e$ is minimal, and isomorphic to im$(e')$ if and only if the idempotents $e$ and $e'$ 
are conjugate in $E$. During this proof we denote the latter situation by $e \sim e'$. 
According to Proposition \ref{prop:7.2} the irreducible quotients of 
$\pi (P,\delta ,\lambda )$ are the modules 
\[
L_{e,\lambda} := L \big( Q, \mr{im}(e), \lambda |_{\mf t^{Q*}} \big) ,
\]
where $e$ runs over the minimal idempotents of $E$. By Theorem \ref{thm:2.1}.b 
$L_{e,\lambda} \cong L_{e',\lambda'}$ if and only if $\lambda = \lambda'$ and $e \sim e'$. 
Let minId$(E)$ be a collection of minimal idempotents of $E$,
containing exactly one element from every conjugacy class. Then
\[
\mr{Irr}_{(P,\delta ,U)} (\mh H ) = \{ L_{e,\lambda} : e \in \mr{minId}(E), \lambda \in U \}
\cong U \times \mr{minId}(E) .
\]
By the second assumption of the theorem
\[
\{ \gamma \in \Gamma : \gamma (P,\delta ,\lambda ) \cong (P,\delta ,\lambda ) \}
\]
does not depend on $\lambda \in U$. The correspoding isomorphism between 
$\pi (P,\delta ,\lambda )$ and $\pi (\gamma (P,\delta ,\lambda ))$ is induced by the map 
\begin{align*}
& J_\gamma : \pi_Q (P,\delta ,\lambda ) \to \pi_Q (\gamma (P,\delta ,\lambda )) \,,\\
& J_\gamma (h \otimes v) = \psi_\gamma (h) \otimes I^\gamma_\delta (v)  
\hspace{2cm} h \in \mh H_Q , v \in V_\delta ,
\end{align*}
with $I^\gamma_\delta$ as in \eqref{eq:8.7}. Hence 
\[
\Gamma_e := \{ \gamma \in \Gamma : L_{e,\lambda} \circ \psi_\gamma^{-1} \cong L_{e,\lambda} \}
= \{ \gamma \in \Gamma : \gamma (P,\delta ,\lambda ) \cong (P,\delta ,\lambda ) ,
J_\gamma e J_\gamma^{-1} \sim e \}
\]
is also independent of $\lambda \in U$. In particular the 2-cocycle 
$\kappa : \Gamma_e \times \Gamma_e \to \mh C^\times$ from \eqref{eq:A.1} measures the lack 
of multiplicativity of $\gamma \mapsto J_\gamma$ and does not depend on $\lambda \in U$.

By Theorem \ref{thm:A.1}.b and Lemma \ref{lem:A.2} the irreducible quotients of 
$\pi' (P,\delta ,\lambda )$ are the $\mh H'$-modules 
$\mr{Ind}_{\Gamma_e \ltimes \mh H}^{\Gamma \ltimes \mh H} (M \otimes L_{e,\lambda})$,
where $e \in \mr{minId}(E)$ and $M$ is an irreducible $\mh C [\Gamma_e ,\kappa]$-module.
So for every $\lambda \in U$ and $e \in \mr{minId}(E)$ we get a packet of $\mh H'$-modules
\[
P_{e,\lambda} := \big\{ \mr{Ind}_{\Gamma_e \ltimes \mh H}^{\Gamma \ltimes \mh H} 
(M \otimes L_{e,\lambda}) : M \in \mr{Irr} \big( \mh C [\Gamma_e ,\kappa ] \big) \big\} .
\]
We note that the ingredients $e$ and $M$ allow us to identify
$\cup_{e \in \mr{minId}(E)} P_{e,\lambda}$ for different $\lambda \in U$. In particular
\begin{equation}\label{eq:10.1}
\bigcup_{\lambda \in U} \bigcup_{e \in \mr{minId}(E)} P_{e,\lambda} \; \cong \;
U \times \bigcup_{e \in \mr{minId}(E)} P_{e,\lambda_0} .
\end{equation}
By Theorem \ref{thm:A.1}.e $P_{e,\lambda}$ and $P_{e',\lambda'}$ contain a common 
$\mh H'$-module if and only if there exists $\gamma \in \Gamma$ such that 
$L_{e,\lambda} \circ \psi_\gamma^{-1} \cong L_{e',\lambda'}$ as $\mh H$-modules. Moreover, 
in this case $P_{e,\lambda} = P_{e',\lambda'}$ as subsets of Irr$(\mh H')$. Therefore 
$\mr{Irr}_{P,\delta ,U} (\mh H')$ consists of the $\Gamma$-equivalence classes in 
\eqref{eq:10.1}.

Notice that $\gamma (\lambda , V_0)$ is only defined if $\gamma (\lambda ) \in U$, 
in which case $\gamma (\lambda , V_0) = (\gamma (\lambda), \gamma (V_0 ))$ for a suitable 
$\gamma (V_0 ) \in \cup_{e \in \mr{minId}(E)} P_{e,\lambda_0 }$. 
Because the action is continuous and $\cup_{e \in \mr{minId}(E)} P_{e,\lambda_0}$ is a finite
set, $\gamma (V_0)$ does not depend on $\lambda$. Since every $\lambda \in U$ has the
same stabilizer in $\Gamma$, we conclude that
\[
\mr{Irr}_{(P,\delta ,U)} (\mh H') \; \cong \; U / \Gamma \times \big( \cup_{e \in \mr{minId}(E)} 
P_{e,\lambda_0} \big) / \Gamma \; \cong \; 
U / W'_\delta \times \mr{Irr}_{(P,\delta ,\lambda_0 )} (\mh H') . \qquad \Box
\]
\vspace{1mm}

The Jacobson topology on $\mr{Irr}_{P,\delta} (\mh H')$ is rather tricky. To put it concisely,
$\mr{Irr}_{P,\delta} (\mh H')$ is a nonseparated affine scheme with a finite-to-one morphism 
onto $\mf t^P / W'_\delta$. Let us also try to describe it more precisely as a space ``fibered'' 
over $\mf t^P / W'_\delta$. Over subsets $U \subset \mf t^P$ as in Proposition \ref{prop:10.1}, 
$\mr{Irr}_{P,\delta} (\mh H')$ looks like $n(U)$ disjoint copies $U_n$ of $U / W'_\delta$. 
If $\lambda \in U$ moves to some boundary point $\mu$ of $U$, 
so that $\mu$ is fixed by more elements than $\lambda$, then $n(\mu )$ may be larger than 
$n(\lambda ) = n(U)$. The topology is such that one path in any of the $U_n$ may converge to 
several points lying over $\mu$. Thus, at boundary of $U$ two different $U_n$'s may either 
converge to the same point, or split up in many disjoint parts, or a combination of both.

In particular over any halfline $\{ r \lambda : r \geq 0 \}$ the space 
$\mr{Irr}_{P,\delta} (\mh H')$ has two typical parts. The fibers over $\{ r \lambda : r > 0 \}$ 
form $n(\lambda )$ disjoint copies of $\mh R_{>0}$, while the fiber over $0 \in \mf t^P$ is 
special. Usually this fiber consists of more than $n(\lambda )$
points, and several of copies of $\mh R_{>0}$ may be connected through the fiber at $0$.

We must also discuss how the different series fit together. By Theorem \ref{thm:8.3} 
$\mr{Irr}_{P,\delta} (\mh H')$ and $\mr{Irr}_{Q,\sigma } (\mh H')$ are either disjoint or 
equivalent. The latter happens if and only if there is some $u \in W'(P,Q)$ with 
$u (\delta ) \cong \sigma$. In view of Theorem \ref{thm:6.2}.b there remain only finitely 
many inequivalent series. Unfortunately a general procedure for constructing discrete series 
representations is lacking, so it is difficult to predict how many series we get.

If $|P| = |Q|$, then the closures of $\mr{Irr}_{P,\delta} (\mh H')$ and $\mr{Irr}_{Q,\sigma} 
(\mh H')$ in the spectrum of $\mh H'$ are either equal or intersect in a subset of smaller 
dimension. One might hope that this subset is necessarily empty, but Theorem \ref{thm:8.3} 
just falls short of proving so.

On the other hand, when $|P| < |Q|$, then the closure of $\mr{Irr}_{P,\delta} (\mh H')$ may 
actually contain $\mr{Irr}_{Q,\sigma } (\mh H')$. This is due to the irreducible constituents 
of $\pi' (P,\delta,\lambda )$ that are not quotients. In fact for every $(Q,\sigma )$ and every 
$P \subsetneq Q$, there exists a $\delta$ for which this happens. In particular the closure 
of $\mr{Irr}_{\emptyset ,\delta_0} (\mh H')$, with $\delta_0$ the irreducible representation 
of $\mh C = \mh H_{\emptyset}$, is the entire spectrum of $\mh H'$.
\vspace{4mm}

\appendix
\section{Appendix: Clifford theory}

Let $A$ be a complex algebra and let $\Gamma$ be a finite group acting on $A$ by algebra 
automorphisms $\psi_\gamma$. Clifford theory describes the relation between the irreducible 
modules of $A$ and those of the crossed product $\Gamma \ltimes A$. For any irreducible 
$A$-module $(\pi ,V_\pi )$ we put
\begin{equation}
\Gamma_\pi = \{ \gamma \in \Gamma : \pi \circ \psi_\gamma^{-1} \cong \pi \} \,.
\end{equation}
For every $\gamma \in \Gamma_\pi$ we choose an isomorphism $I^\gamma$ between $\pi$ and 
$\pi \circ \psi_\gamma^{-1}$:
\begin{align*}
& I^\gamma : V_\pi \to V_\pi , \\
& I^\gamma (\pi (h ) v) = \pi (\psi_\gamma^{-1} h) I^\gamma (v) \qquad h \in A , v \in V_\pi \,.
\end{align*}
We assume that $V_\pi$ has (at most) countable dimension, so that Schur's lemma assures us that 
$I^\gamma$ is unique up to a nonzero complex number. This gives rise to a 2-cocycle
\begin{equation}\label{eq:A.1}
\begin{split}
& \kappa : \Gamma_\pi \times \Gamma_\pi \to \mh C^\times , \\
& I^{\gamma \gamma'} = \kappa (\gamma ,\gamma' ) I^\gamma I^{\gamma'} .
\end{split}
\end{equation}
Let $\mh C [\Gamma_\pi ,\kappa]$ be the associative $\mh C$-algebra with basis 
$\{ T_\gamma : \gamma \in \Gamma_\pi \}$ and multiplication defined by
\[
T_\gamma T_{\gamma'} = \kappa (\gamma ,\gamma') T_{\gamma \gamma'} \,.
\]
Different choices of $I^\gamma$ give rise to different $\kappa$'s, but to isomorphic 
twisted group algebras $\mh C [\Gamma_\pi ,\kappa]$. Hence we may and will assume that 
$I^e = \mr{Id}_{V_\pi}$, so that $T_e$ is the unit of this algebra. Now 
$\mh C [\Gamma_\pi ,\kappa ] \otimes V_\pi$ becomes a $\Gamma_\pi \ltimes A$-module by
\[
(\gamma \otimes h) (k \otimes v) = T_\gamma k \otimes I^\gamma (\pi (h) v) \qquad
\gamma \in \Gamma_\pi , h \in A , k \in \mh C [\Gamma_\pi ,\kappa ] , v \in V_\pi .
\]

\begin{thm}\label{thm:A.1} \textup{(Clifford theory)}
\begin{description}
\item[a)] There is an isomorphism of $\Gamma_\pi \ltimes A$-modules
\begin{align*}
& T : \mr{Ind}_{A}^{\Gamma_\pi \ltimes A} V_\pi \to 
\mh C [\Gamma_\pi ,\kappa ] \otimes V_\pi , \\
& T (\gamma \otimes v ) = T_\gamma \otimes I^\gamma (v) .
\end{align*}
\item[b)] The map $M \mapsto \mr{Ind}_{\Gamma_\pi \ltimes A}^{\Gamma \ltimes A} 
(T^{-1} (M \otimes V_\pi ))$ is an isomorphism between the following categories:
\begin{itemize}
\item submodules of the left regular representation of $\mh C [\Gamma_\pi ,\kappa ]$,
\item $\Gamma \ltimes A$-submodules of $\mr{Ind}_{A}^{\Gamma \ltimes A} V_\pi$.
\end{itemize}
\item[c)] The $\Gamma \ltimes A$-module $\mr{Ind}_{A}^{\Gamma \ltimes A} V_\pi$
is completely reducible.
\item[d)] Every irreducible $\Gamma \ltimes A$-module can be obtained from this construction.
\item[e)] Let $(\pi ,V_\pi )$ and $(\rho ,V_\rho )$ be irreducible $A$-modules whose 
inductions to $\Gamma \ltimes A$ have a common irreducible summand. Then 
\[
\mr{Ind}_{A}^{\Gamma \ltimes A} V_\pi \cong 
\mr{Ind}_{A}^{\Gamma \ltimes A} V_\rho
\]
as $\Gamma \ltimes A$-modules, and there exists $\gamma \in \Gamma$ with 
$\pi \circ \psi_\gamma^{-1} \cong \rho$.
\end{description}
\end{thm}
\emph{Proof.}
This theorem is by no means original, similar results can be found for example in 
\cite{Mac} and \cite[p. 24]{RaRa}.\\
a) is a simple direct verification.\\
b) Clearly the given map is functorial. To construct its inverse, consider
\begin{equation}\label{eq:A.7}
\mr{Ind}_{\Gamma_\pi \ltimes A}^{\Gamma \ltimes A} (\mh C [\Gamma_\pi ,\kappa ] \otimes V_\pi )
= {\ts \bigoplus_{g \in \Gamma / \Gamma_\pi}} g (\mh C [\Gamma_\pi ,\kappa ] \otimes V_\pi )
\end{equation}
as an $A$-module. By the definition of $\Gamma_\pi$, the right hand side is precisely
the decomposition into isotypical components. So for any $A$-submodule $N$ of \eqref{eq:A.7}:
\[
N = {\ts \bigoplus_{g \in \Gamma / \Gamma_\pi}} 
N \cap g (\mh C [\Gamma_\pi ,\kappa ] \otimes V_\pi ) .
\]
If moreover $N$ is a $\Gamma \ltimes A$-submodule, then 
\[
N = {\ts \bigoplus_{g \in \Gamma / \Gamma_\pi}} 
    g (N \cap \mh C [\Gamma_\pi ,\kappa ] \otimes V_\pi )
= \mr{Ind}_{\Gamma_\pi \ltimes A}^{\Gamma \ltimes A} 
(N \cap \mh C [\Gamma_\pi ,\kappa ] \otimes V_\pi ) .
\]
Since it is an $A$-submodule,
\[
N \cap (\mh C [\Gamma_\pi ,\kappa ] \otimes V_\pi ) =  
M \otimes V_\pi \quad \subset \quad \mh C [\Gamma_\pi ,\kappa ] \otimes V_\pi ,
\]
for some subspace $M \subset \mh C [\Gamma_\pi ,\kappa ]$. Then $M \otimes V_\pi$ is stable 
under $\Gamma_\pi$, so $M$ has to be a $\mh C [\Gamma_\pi ,\kappa ]$-submodule. \\
c) $\mh C [\Gamma_\pi ,\kappa ]$ is semisimple, so its left regular representation is 
completely reducible. Combine this with b). \\
d) Let $N$ be any irreducible $\Gamma \ltimes A$-module, and let $V_\pi$ be an irreducible 
$A$-submodule of $N$. By Frobenius reciprocity 
\[
\mr{Hom}_{\Gamma \ltimes A} \big( \mr{Ind}_{A}^{\Gamma \ltimes A} V_\pi , N \big)
\cong \mr{Hom}_{A} (V_\pi ,N) \neq 0 ,
\]
so $N$ is isomorphic to a quotient of the $\Gamma \ltimes A$-module 
$\mr{Ind}_{A}^{\Gamma \ltimes A} V_\pi$. By c) $N$ is also isomorphic to a direct summand
of $\mr{Ind}_{A}^{\Gamma \ltimes A} V_\pi$.\\
e) As $A$-modules
\begin{equation}\label{eq:A.10}
\begin{split}
& \mr{Ind}_{A}^{\Gamma \ltimes A} V_\pi \cong {\ts \bigoplus_{\gamma \in \Gamma}} \pi \circ 
\psi_\gamma^{-1} , \\
& \mr{Ind}_{A}^{\Gamma \ltimes A} V_\rho \, \cong {\ts \bigoplus_{\gamma \in \Gamma}} \rho \circ 
\psi_\gamma^{-1} .
\end{split}
\end{equation}
By assumption these modules have a common $A$-submodule, which implies that they are actually
isomorphic. In particular, for some $\gamma \in \Gamma$ there exists a linear bijection
$\phi : V_\pi \to V_\rho$ that provides an isomorphism between $\pi \circ \psi_\gamma^{-1}$ 
and $\rho$. Then 
\begin{equation}\label{eq:A.8}
\mr{Ind}_{A}^{\Gamma \ltimes A} V_\pi \to \mr{Ind}_{A}^{\Gamma \ltimes A} V_\rho 
: g \otimes v \mapsto g \gamma^{-1} \otimes \phi (v)
\end{equation}
is an isomorphism of $\Gamma \ltimes A$-modules. $\qquad \Box$
\\[2mm]

As a direct consequence of part b), the number of inequivalent irreducible constituents of the
$\Gamma \ltimes A$-module $\mr{Ind}_{A}^{\Gamma \ltimes A} V_\pi$ depends only on two things:
the group $\Gamma_\pi$ and the cocycle $\kappa$.

\begin{lem}\label{lem:A.2}
Let $V$ be any $A$-module. The irreducible quotients of $\mr{Ind}_{A}^{\Gamma \ltimes A} V$ 
are precisely the irreducible summands of the modules $\mr{Ind}_{A}^{\Gamma \ltimes A} Q$, 
where $Q$ runs over the irreducible quotients of $V$.
\end{lem}
\emph{Proof.}
Let $V' = \mr{Ind}_{A}^{\Gamma \ltimes A} (V) / N$ be an irreducible quotient $\Gamma 
\ltimes A$-module. Any proper $A$-submodule of $Q := V / (N \cap V)$ would generate a proper 
$\Gamma \ltimes A$-submodule of $V'$, which by assumption is impossible. Therefore $Q$ is an 
irreducible quotient of $V$. By construction $V'$ is a quotient $\Gamma \ltimes A$-module of 
\[
\mr{Ind}_{A}^{\Gamma \ltimes A} Q = 
\bigoplus\nolimits_{\gamma \in \Gamma} \gamma V / (N \cap \gamma V) ,
\]
and by Theorem \ref{thm:A.1}.b it is actually a summand. $\qquad \Box$

\end{document}